\documentclass[leqno,CJK]{siamltex704}
\usepackage{amssymb,amsmath,graphicx,amscd,mathrsfs}
\usepackage{color,xcolor,amsmath}
\usepackage{amsmath}
\usepackage{graphicx}
\usepackage{mathrsfs}
\usepackage{float}
\usepackage{amsfonts,amssymb}
\usepackage{dsfont}
\usepackage{pifont}
\usepackage{hyperref}
\usepackage{multirow}
\usepackage{placeins}
\usepackage{bm}
\usepackage{enumerate}
\usepackage{multirow}
\usepackage{multicol}
\usepackage{subcaption}
\numberwithin{equation}{section}
\def\3bar{{|\hspace{-.02in}|\hspace{-.02in}|}}

\def\b0{\boldsymbol{0}}

\def\bf{{\mathbf{f}}}

\newtheorem{example}{Example}[section]
\newtheorem{remark}{Remark}[section]

\newtheorem{assumption}{Assumption}

\allowdisplaybreaks 

\begin{document}

\title{Spectral properties of high dimensional rescaled sample correlation matrices}

\author{
Weijiang Chen\thanks{Jilin University, Changchun,
China (chenwj21@mails.jlu.edu.cn).}
\and
Shurong Zheng\thanks{Northeast Normal University, Changchun, China 
 (zhengsr@nenu.edu.cn).The research of this author was supported in part by the National
Natural Science Foundation of China (NSFC 12231011, 12326606).}
\and
Tingting Zou\thanks{Jilin University, Changchun, China 
 (zoutt260@jlu.edu.cn).
}
}
\maketitle
\begin{abstract}
High-dimensional sample correlation matrices are a crucial class of random matrices in multivariate statistical analysis. The central limit theorem (CLT) provides a theoretical foundation for statistical inference. In this paper, assuming that the data dimension increases proportionally with the sample size, we derive the limiting spectral distribution of the matrix $\widehat{\mathbf{R}}_n\mathbf{M}$ and establish the CLTs for the linear spectral statistics (LSS) of $\widehat{\mathbf{R}}_n\mathbf{M}$ in two structures: linear independent component structure and elliptical structure. In contrast to existing literature, our proposed spectral properties do not require $\mathbf{M}$ to be an identity matrix. Moreover, we also derive the joint limiting distribution of LSSs of $\widehat{\mathbf{R}}_n \mathbf{M}_1,\ldots,\widehat{\mathbf{R}}_n \mathbf{M}_K$.
As an illustration, an application is given for the CLT.
\end{abstract}
\begin{keywords}
random matrix, multivariate statistical analysis, central limit theorem, linear independent component structure, elliptical structure, joint limiting distribution.
\end{keywords}


\section{Introduction}
\numberwithin{equation}{section}
Sample correlation matrix is an important matrix in multivariate statistical analysis. Suppose that $\mathbf{y}_{1},\ldots, \mathbf{y}_{n}$ are \emph{independent and identically distributed} (i.i.d.) from a \emph{p}-dimensional population $\mathbf{y}$ with a mean vector $\boldsymbol{\mu}$ and a covariance matrix $\boldsymbol{\Sigma}$. The population correlation matrix is defined as
\begin{equation*}
    \mathbf{R}=\left[\textup{diag}(\boldsymbol{\Sigma})\right]^{-1/2}\boldsymbol{\Sigma}\left[\textup{diag}(\boldsymbol{\Sigma})\right]^{-1/2},
\end{equation*}
where $\textup{diag}(\boldsymbol{\Sigma})$ is a diagonal matrix consisting of the diagonal elements of $\boldsymbol{\Sigma}$. The sample covariance matrix $\mathbf{S}_{n}$ and sample correlation matrix $\widehat{\mathbf{R}}_{n}$ are defined as
\begin{equation}
    \mathbf{S}_{n}=(n-1)^{-1}\displaystyle\sum_{j=1}^{n}(\mathbf{y}_{j}-\overline{\mathbf{y}})(\mathbf{y}_{j}-\overline{\mathbf{y}})^{\top}, ~~\widehat{\mathbf{R}}_{n}=\left[\textup{diag}( \mathbf{S}_{n})\right]^{-1/2}\mathbf{S}_{n}\left[\textup{diag}( \mathbf{S}_{n})\right]^{-1/2},
\end{equation}
where $\overline{\mathbf{y}}=n^{-1}\sum_{j=1}^{n}\mathbf{y}_{j}$ is the sample mean.

\cite{fan2022estimating} showed that it is important and necessary to study the random matrix theory of sample correlation matrices. Many studies have examined the spectral properties of high-dimensional sample matrices $\widehat{\mathbf{R}}_n$. For the \emph{limiting spectral distribution} (LSD) and extreme eigenvalues of $\widehat{\mathbf{R}}_n$, \cite{jiang2004limiting} derived the M-P law and almost sure convergence of the largest eigenvalue under $\mathbf{R}=\mathbf{I}_p$.
Subsequently, \cite{xiao2010almost} verified the almost sure convergence of the smallest eigenvalue. \cite{bao2012tracy} and \cite{pillai2012edge} simultaneously showed that the extreme eigenvalues converge in distribution to Tracy-Widom law. \cite{heinyandyao2022} found the LSD in the heavy-tailed case that the sample has infinite variace. For a general $\mathbf{R}$, \cite{el2009concentration} and \cite{yin2023central} derived the LSD of $\widehat{\mathbf{R}}_n$ under linear independent component and elliptical structures, respectively.

As for the \emph{central limit theorem} (CLT) for \emph{linear spectral statistics} (LSS) of a high-dimensional sample correlation matrix $\widehat{\mathbf{R}}_{n}$, \cite{gao2017high} established the CLT for the LSS in the case where $\mathbf{R}=\mathbf{I}_{p}$. \cite{mestre2017correlation} relaxed the restriction $\mathbf{R}=\mathbf{I}_{p}$ but required a Gaussian assumption. \cite{jiang2019determinant} derived a CLT for the logarithm of the determinant of $\widehat{\mathbf{R}}_n$ under Gaussian population. \cite{yin2022spectral} established the CLT for the LSS of the rescaled sample correlation matrix $\widehat{\mathbf{R}}_{n}\mathbf{R}^{-1}$ when the sample has an independent component structure. 
\cite{yin2023central} obtained a general CLT for the LSS of $\widehat{\mathbf{R}}_{n}$ under both the independent component structure and elliptical structure assumptions. \cite{yin2022properties} derived the central limit theorem for the linear statistics of the eigenvectors of $\widehat{\mathbf{R}}_n$ under a general fourth moment condition.

In order to test independence of $p$-variates of the population under high dimensional settings, a number of studies have been conducted. \cite{jiang2004asymptotic}, \cite{zhou2007asymptotic}, \cite{liu2008asymptotic}, \cite{cai2011limiting}, and \cite{cai2012phase} derived the asymptotic distribution of maximum-norm-type statistics based on the largest entries of $\widehat{\mathbf{R}}_{n}$. \cite{schott2005testing} and \cite{mestre2017correlation} established the asymptotic behavior of Frobenius-norm-type statistic $\textup{tr}[(\widehat{\mathbf{R}}_{n}-\mathbf{I}_p)^2]$. To cope with heavy-tailed population, \cite{bao2015spearman} constructed test methods based on the polynomial functions of the spectrum of Spearman's rank correlation matrices. \cite{han2017rank} proposed two families of maximum-norm-type rank statistics including Spearman's rho and Kendall's tau as special cases. \cite{leung2018testing} considered three types of test statistics consisting of sums or sums of squares of pairwise rank correlations. \cite{li2021central} established the CLT for the LSS of Kendall's rank correlation matrices and applied the new CLT to construct test methods.

For testing whether the population correlation matrix equals to a given matrix, that is, testing the hypothesis 
\begin{equation}\label{test}
    H_0: \mathbf{R}=\mathbf{R}_0~~~v.s.~~~H_1:\mathbf{R}\neq\mathbf{R}_0,
\end{equation}
where $\mathbf{R}_0$ is a pre-specified matrix, the relevant literature is relatively rare. \cite{zheng2019test} proposed the test statistic $T_2=\textup{tr}[(\widehat{\mathbf{R}}_{n}-\mathbf{R}_0)^{2}]$ and derived its the asymptotic distributions under both null and alternative hypotheses. \cite{yin2022spectral} used $T_1=\textup{tr}[\widehat{\mathbf{R}}_{n}\mathbf{R}_0^{-1}-\mathbf{I}_p]^{2}$ to construct the testing method. In fact, both these two methods have their own advantages and their relative performance varies case by case, which will be shown in the subsequent simulation. 

Consider the following two scenarios:
\begin{itemize}
    \item Scenario~$1$:~$\boldsymbol{\Sigma}=\boldsymbol{\Gamma}\boldsymbol{\Gamma}^{\top}$, where $\boldsymbol{\Gamma}=\mathbf{I}_p+\theta \mathbf{A}$, and all elements of $\mathbf{A}=(a_{ij})_{i,j=1,\ldots,p}$ are independent and identically distributed (i.i.d.) from the uniform distribution $U(-p^{-2/3},p^{-2/3})$;
    \item Scenario~$2$:~ $$\boldsymbol{\Sigma}=\mathbf{U }(\mathbf{I}_p+\mathbf{D})\mathbf{U}^{\top}+\theta \mathbf{1}_p^{\top}\mathbf{1}_p,$$ 
    where $\mathbf{U}$ is the eigenvector matrix of $\mathbf{Z}^{\top}\mathbf{Z}$ with the 
    elements $z_{ij}(i,j=1,\ldots,p)$ being i.i.d. from $N(0,1)$, $\mathbf{D}=\textup{diag}(d_{11},d_{22},\ldots,d_{pp})$ is a diagonal matrix, and the elements of $\mathbf{D}$ are i.i.d. from the uniform distribution $U(0,1)$.
\end{itemize}
The parameter setting is as follows:
\begin{itemize}
    \item Dimension: $p=100, 150, 200, 250, 300, 350, 400;$
    \item Ratio of dimension and sample size: $y_n= p/n =0.5;$
    \item We set $\theta=0$ under the null hypothesis $H_0$ in both models, $\theta=1$ and $\theta=0.1$ to evaluate empirical powers in Scenario $1$ and Scenario $2$, respectively. 
\end{itemize}
For each setting, we run simulation $10000$ times and calculate the corresponding sample means $\mu_i$ and sample variances $\sigma_i$ for $i=1,2$, the rejection regions of the statistics $T_1, T_2$ for testing (\ref{test}) at the level $5\%$ are as follows:
\begin{align*}
 &\{\mathbf{y}_1,\ldots,\mathbf{y}_n:\sigma_1^{-1}|T_1-\mu_1|>q_{0.975}\},\\
 &\{\mathbf{y}_1,\ldots,\mathbf{y}_n:\sigma_2^{-1}|T_1-\mu_2|>q_{0.975}\},
\end{align*}
where $q_{0.975}$ is the $97.5\%$ quantile of $N(0,1)$.

\begin{figure}[h]
        \centering
        \includegraphics[width=0.6\linewidth]{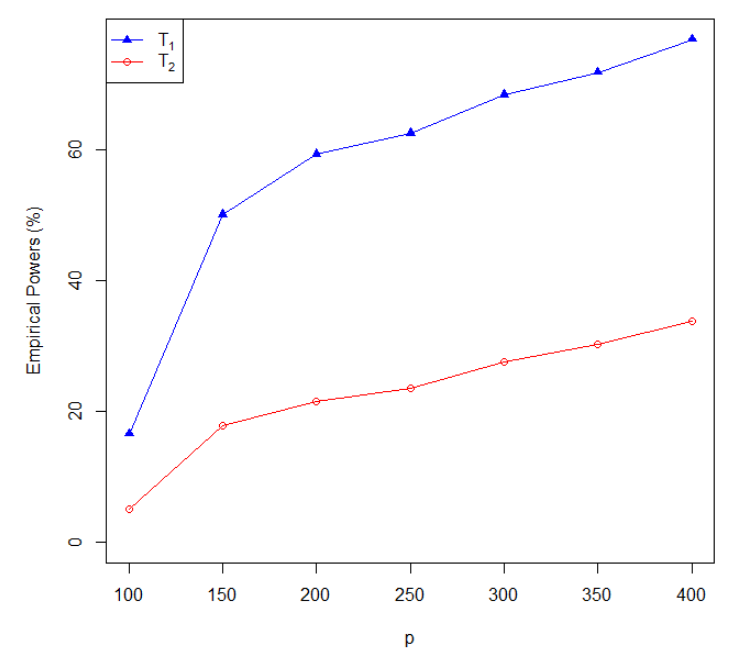}
        \caption{Empirical powers of $T_1$ and $T_2$ in Scenario $1$. The ratio of dimension and sample size is 0.5.}
        \label{fig:T1good}
\end{figure}
\begin{figure}[h]
        \centering
        \includegraphics[width=0.6\linewidth]{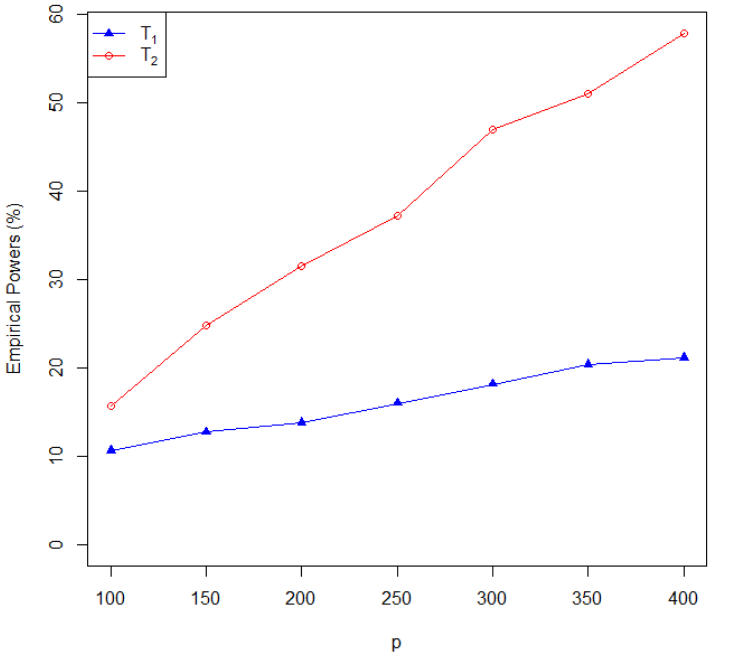}
        \caption{Empirical powers of $T_1$ and $T_2$ in Scenario $2$. The ratio of dimension and sample size is 0.5.}
        \label{fig:T2good}
\end{figure}

Fig.1 shows the empirical powers of two tests in Scenario $1$. In this setting, $T_1$ is more powerful than $T_2$; Fig.2 displays the performances of two tests in Scenario $2$. In this setting, $T_2$ performs better than $T_1$. In order to combine the advantages of $T_1$ and $T_2$, we consider the test statistic as follows
$$
T_m=\max\left\{\sigma_1^{-1}|T_1-\mu_1|,~~\sigma_2^{-1}|T_2-\mu_2|\right\}.
$$
Note that determining the asymptotic distribution of $T_m$ depends on the joint limiting distribution of the LSSs of $\widehat{\mathbf{R}}_n\mathbf{R}^{-1}_0$, $\widehat{\mathbf{R}}_n$, and $\widehat{\mathbf{R}}_n\mathbf{R}_0$. Therefore, the main task in this paper is to establish the central limit theorem for LSSs of $\widehat{\mathbf{R}}_n\mathbf{M}_1,\ldots,\widehat{\mathbf{R}}_n\mathbf{M}_K$, where $\mathbf{M}_1,\ldots,\mathbf{M}_K$ are pre-specified matrices.

The arrangement of this paper is as follows: Section 2 derives the limiting spectral distributions of the rescaled sample correlation matrix $\widehat{\mathbf{R}}_{n}\mathbf{M}$. Section 3 establishes the central limit theorem of LSSs for the rescaled sample correlation matrix $\widehat{\mathbf{R}}_{n}\mathbf{M}$ under the elliptical structure. Section 4 establishes the central limit theorem of LSSs for $\widehat{\mathbf{R}}_{n}\mathbf{M}$ under the independent component structure. The joint limiting distribution of LSSs for $\widehat{\mathbf{R}}_{n}\mathbf{M}_1,\ldots,\widehat{\mathbf{R}}_{n}\mathbf{M}_K$ is given in Section 5. Section 6 provides an application to illustrate the CLT of the LSS for $\widehat{\mathbf{R}}_{n}\mathbf{M}
_1,\ldots,\widehat{\mathbf{R}}_{n}\mathbf{M}_K$.

\section{Limiting spectral distribution of $\widehat{\mathbf{R}}_{n}\mathbf{M}$}
To derive the limiting spectral distribution of $\widehat{\mathbf{R}}_{n}\mathbf{M}$, we need some assumptions.
\numberwithin{equation}{section}
\begin{assumption}
    Assume that the i.i.d. samples $\mathbf{y}_{1},\ldots,\mathbf{y}_{n}$ satisfy the following elliptical structure:
    \begin{equation}\label{eq2.1}
       \mathbf{y}_{j}=\rho_{j} \boldsymbol{\Gamma}\mathbf{x}_{j}+\bm{\mu},~j=1,\ldots,n,
    \end{equation}
where $\boldsymbol{\Gamma}$ is a $p\times p $ non-random  matrix with rank$(\boldsymbol{\Gamma})=p $, $\mathbf{x}_{j}$ is a $p$-dimensional random direction independent of $\rho_{j}$~and uniformly distributed on the unit sphere $S^{p-1}$ in $\mathbb{R}^{p}$, and $\rho_{j}$ is a non-negative random radius satisfying
\begin{equation}
    \textup{E}\rho_{j}^{2}=p,~\textup{E}\rho_{j}^{4}=p^{2}+\tau p +o(p),~\textup{E}\left\vert\dfrac{\rho_{j}^{2}-p}{\sqrt{p}} \right\vert ^{2+\varepsilon} < \infty,
\end{equation}
for constants $\tau \geq 0$ and $\varepsilon>0.$
\end{assumption}
\begin{assumption}
Assume that the i.i.d. samples $\mathbf{y}_{1},\ldots,\mathbf{y}_{n}$ satisfy the following linear independent component structure:
\begin{equation}
       \mathbf{y}_{j}=\boldsymbol{\Gamma}\mathbf{x}_{j}+\boldsymbol{\mu},~j=1,\ldots,n,
    \end{equation}
where $\boldsymbol{\Gamma}$ is a $p\times p $ non-random  matrix, $\textup{E}\mathbf{y}_{j}=\boldsymbol{\mu},$ and $\mathbf{x}_{j}=(x_{1j},\ldots, x_{pj})^{\top}$ is a $p$-dimensional random vector with i.i.d. entries satisfying
\begin{equation}
    \textup{E}x_{ij}=0,~\textup{E}x_{ij}^{2}=1,~\textup{E}x_{ij}^{4}=\beta_{x}+3+o(1), \textup{E}(\vert x_{ij}\vert^{4}(\textup{log}|x_{ij}|)^{2+2\varepsilon})<\infty.
\end{equation}
\end{assumption}
\begin{assumption}
Let $\mathbf{G}=[\textup{diag}(\boldsymbol{\Sigma})]^{-\frac{1}{2}}\boldsymbol{\Gamma}.$ Assume that the empirical spectral distribution ({\sl ESD}) $H_n=p^{-1}\sum_{i=1}^{p}\delta(\lambda^{{\mathbf{RM}}}_{i}\leq x)$ of $\mathbf{RM}=\mathbf{GG^{\top}M}$ weakly converges to a proper distribution $H$ as $p\to \infty$
, where $\delta$ is an indicator function and $\lambda^{{\mathbf{RM}}}_i$ is the $i$th largest eigenvalue of ${\mathbf{RM}}$. Moreover, the spectral norms of $\mathbf{R}, \mathbf{M}, \mathbf{M^{-1}}$ are uniformly bounded in $p$.
\end{assumption}
\begin{assumption}
A convergence regime is required as $y_{n}=p/n\to y\in (0,+\infty).$
\end{assumption}
\begin{assumption}
The functions $g_{1},\ldots, g_{K}$ are known analytic functions in a domain containing
    \begin{equation*}
        \left[\displaystyle\liminf_{p}\lambda_{\min}^{\mathbf{RM}}\cdot \delta_{\{0\leq y\leq 1\}}(1-\sqrt{y})^{2}, \displaystyle\limsup_{p}\lambda_{\max}^{\mathbf{RM}}\cdot (1+\sqrt{y})^{2}\right]
    \end{equation*}
    where $\lambda_{\min}^{\mathbf{RM}}$ and $\lambda_{\max}^{\mathbf{RM}}$ are the minimum and maximum eigenvalues of $\mathbf{RM}$, respectively.
\end{assumption}

Assumption $3$ requires the bounded spectral norm of non-random matrices $\mathbf{R}$, $\mathbf{M}$ and $\mathbf{M}^{-1}$. Assumption $4$ gives the convergence regime of the dimension $p$ and sample size $n$.

The following theorem provides the LSD $F^{y,H}(x)$ of $F_{n}(x)$.
\begin{theorem}\label{th2.1}
    Under Assumptions $1-3-4$ or $2-3-4$, the ESD $F_{n}(x)$ of
 $\widehat{\mathbf{R}}_{n}\mathbf{M}$ converges almost surely to the LSD $F^{y,H}$, whose Stieltjes transform $s(z)$ is the only solution to the equation
    \begin{equation}\label{eq2.5}
        s(z)=\displaystyle\int\dfrac{1}{t\left[1-y-yzs(z)\right]-z}\mathrm{d}H(t),~~z\in \mathbb{C}^{+},
    \end{equation}
    in the set $\{ s(z): -(1-y)/z+ys(z)\in \mathbb{C}^{+}\}$, where $\mathbb{C}^{+}=\{z\in \mathbb{C}: \mathfrak{I}(z)>0\}$ with $\mathfrak{I}(z)$ being the imaginary part of $z$. Letting
    \begin{equation*}
        \underline{s}(z)=-\dfrac{1-y}{z}+ys(z),~~z\in \mathbb{C}^{+},
    \end{equation*}
    then (\ref{eq2.5}) can be re-expressed as
    \begin{equation*}
        z=-\dfrac{1}{\underline{s}(z)}+y\displaystyle\int \dfrac{t}{1+t\underline{s}(z)}\mathrm{d}H(t).
    \end{equation*}
\end{theorem}
Let $[a,b]$ be the support set of the LSD $F^{y,H}(x)$; then, $F^{y,H}(x)=0, x<0 $ and
\begin{equation*}
 F^{y,H}(x)= \left\{
\begin{aligned}
 &\displaystyle\int_{0}^{x}f^{y,H}(t)\textup{d}t,~~~~~~~~~~~~~~~~~~~~~~~~~~~~~~~~y\leq 1, x\geq 0,\\
  & \displaystyle\int_{0}^{x}f^{y,H}(t)\textup{d}t+(1-1/y)\delta_{\{x\geq 0\}},~~~~~~~~y>1, x\geq 0,
\end{aligned}
\right.
\end{equation*}
with the limiting spectral density being
\begin{equation}\label{eq2.6}
    f^{y,H}(x)=(y\pi)^{-1}\displaystyle\lim_{z\to x} \mathfrak{I}(\underline{s}(z))\delta_{\{0\leq a\leq x\leq b\}}.
\end{equation}
Note that the Stieltjes transform $s(z)$ of the LSD of  $\widehat{\mathbf{R}}_{n}\mathbf{M}$ are same for elliptical structure and independent component structure.


\section{{\color{black}Central limit theorem of linear spectral statistics of $\widehat{\mathbf{R}}_{n}\mathbf{M}$ under elliptical structure}}

Define the LSS of $\widehat{\mathbf{R}}_{n}\mathbf{M}$ as
\begin{equation}
L_{g_{\ell}}=\displaystyle\sum_{i=1}^{p}g_{\ell}(\widehat{\lambda}_{i}),~~\ell=1,\ldots, K,
\end{equation}
where $g_{\ell}(\cdot), \ell=1,\ldots, K $ are some known analytic functions. We denote $M_{n}(z)=p(s_n(z)-s_{y_{n}}(z)),$ where $s_n(z)$ is the Stieltjes transform of the ESD of $\widehat{\mathbf{R}}_n\mathbf{M}$ and $s_{y_{n}}(z)$ is the Stieltjes transform of the distribution $F^{y_n,H_n}$. To establish the CLT of the LSS of $\widehat{\mathbf{R}}_{n}\mathbf{M}$, for a fixed $K$ and known functions $g_{1},\ldots, g_{K}$, we consider the $K$-dimensional random vector $(W(g_{1}),\ldots, W(g_{K}))$, where
\begin{equation}
W(g_{\ell})=\displaystyle\sum_{i=1}^{p}g_{\ell}(\widehat{\lambda}_{i})\delta_{\{\widehat{\lambda}_{i}>0\}}-p\displaystyle\int_{a}^{b} g_{\ell}(x)f^{y_{n-1},H_n}(x)\mathrm{d}x,~~ \ell=1,\ldots, K 
\end{equation}
and $f^{y_{n-1},H_n}$ is defined in (\ref{eq2.6}) with $y_{n-1}=p/(n-1)$.

This section establishes the CLT of the random vector $(W(g_{1}),\ldots, W(g_{K}))$ under the elliptical structure assumption $1$.
\begin{theorem}\label{th3.1}
Under Assumption $1$ and Assumptions $3-5$ , the random vector $(W(g_{1}),\ldots\\, W(g_{K}))$ weakly converges to a multivariate Gaussian vector $(X_{g_{1}},\ldots, X_{g_{K}})$ with the mean and covariance functions as follows:
\begin{equation}\label{temp1}
\textup{E}X_{g_{\ell}}=-\dfrac{1}{2\pi i}\displaystyle\oint_{\mathcal{C}}g_{\ell}(z)\textup{E}M(z)\mathrm{d}z,
\end{equation}
\begin{equation}\label{temp2}
\textup{Cov}(X_{g_{\ell_{1}}},X_{g_{\ell_{2}}})=-\dfrac{1}{4\pi^{2}}\displaystyle\oint_{\mathcal{C}_{1}}\oint_{\mathcal{C}_{2}}g_{\ell_{1}}(z_{1})g_{\ell_{2}}(z_{2})\textup{Cov}(M(z_{1}),M(z_{2}))\mathrm{d}z_{2}\mathrm{d}z_{1},
\end{equation}
where
\begin{align}\label{eq3.3}
\textup{E}M(z)=&y\displaystyle\int\dfrac{[t\underline{s}'(z)]^{2}}{\underline{s}(z)[1+t\underline{s}(z)]^{3}}\mathrm{d}H(t)+(\tau -2)[1+z\underline{s}(z)]\displaystyle\int\dfrac{t\underline{s}'(z)}{[1+t\underline{s}(z)]^{2}}\mathrm{d}H(t)\\
&+\displaystyle\lim_{n\to \infty } \dfrac{1}{n} \sum_{k=1}^{p}\dfrac{\partial}{\partial z}\left[\underline{s}(z)\mathbf{e}_{k}^{\top}\mathbb{R}(z)\mathbf{R}\mathbf{e}_{k}\cdot \mathbf{e}_{k}^{\top}\mathbf{RM}\mathbb{R}(z)\mathbf{e}_{k}\right]\nonumber\\
&+\displaystyle\lim_{n\to \infty } \dfrac{1}{n} \sum_{k=1}^{p}\dfrac{\partial}{\partial z}\left[\underline{s}(z)\mathbf{e}_{k}^{\top}\mathbf{M}\mathbb{R}(z)\mathbf{R}\mathbf{e}_{k}\cdot \mathbf{e}_{k}^{\top}\mathbf{RM}\mathbb{R}(z)\mathbf{M}\mathbf{e}_{k}\right]\nonumber\\
&+\dfrac{1}{4z}\displaystyle\lim_{n\to \infty}\dfrac{1}{n}\sum_{\ell=1}^{p}\left[\mathbf{e}_{\ell}^{\top}\mathbb{R}(z)\mathbf{e}_{\ell}\right]
+\dfrac{1}{4z}\displaystyle\lim_{n\to \infty}\dfrac{1}{n}\sum_{\ell=1}^{p}\left[\mathbf{e}_{\ell}^{\top}\mathbf{M}\mathbb{R}(z)\mathbf{M}^{-1}\mathbf{e}_{\ell}\right]\nonumber\\
&-\dfrac{1}{4z}\displaystyle\lim_{n\to \infty}\dfrac{1}{n}\sum_{\ell=1}^{p}\left[\mathbf{e}_{\ell}^{\top}\mathbb{R}^{2}(z)\mathbf{e}_{\ell}\right]
-\dfrac{1}{4z}\displaystyle\lim_{n\to \infty}\dfrac{1}{n}\sum_{\ell=1}^{p}\left[\mathbf{e}_{\ell}^{\top}\mathbf{M}\mathbb{R}^{2}(z)\mathbf{M}^{-1}\mathbf{e}_{\ell}\right]\nonumber\\
&-\dfrac{1}{2z}\displaystyle\lim_{n\to \infty}\dfrac{1}{n}\sum_{k,\ell=1}^{p}r_{k\ell}^{2}\mathbf{e}_{\ell}^{\top}\mathbf{M}^{-1}\mathbf{e}_{k}\left[\mathbf{e}_{k}^{\top}\mathbf{M}\mathbb{R}(z)\mathbf{e}_{\ell}\right]\nonumber\\\
&+\dfrac{1}{2z}\displaystyle\lim_{n\to \infty}\dfrac{1}{n}\sum_{k,\ell=1}^{p}r_{k\ell}^{2}\mathbf{e}_{\ell}^{\top}\mathbf{M}^{-1}\mathbf{e}_{k}\left[\mathbf{e}_{k}^{\top}\mathbf{M}\mathbb{R}^{2}(z)\mathbf{e}_{\ell}\right]\nonumber\\
&+\dfrac{1}{4}\displaystyle\lim_{n\to \infty}\dfrac{1}{n}\sum_{k,\ell=1}^{p}r_{k\ell}^{2}\dfrac{\partial}{\partial z}\left[\mathbf{e}_{k}^{\top}\mathbb{R}(z)\mathbf{e}_{\ell}\cdot \mathbf{e}_{\ell}^{\top}\mathbb{R}(z)\mathbf{e}_{k}\right]\nonumber\\
&+\dfrac{1}{2}\displaystyle\lim_{n\to \infty}\dfrac{1}{n}\sum_{k,\ell=1}^{p}r_{k\ell}^{2}\dfrac{\partial}{\partial z}\left[\mathbf{e}_{k}^{\top}\mathbb{R}(z)\mathbf{M}^{-1}\mathbf{e}_{\ell}\cdot \mathbf{e}_{\ell}^{\top}\mathbf{M}\mathbb{R}(z)\mathbf{e}_{k}\right]\nonumber\\
&+\dfrac{1}{4}\displaystyle\lim_{n\to \infty}\dfrac{1}{n}\sum_{k,\ell=1}^{p}r_{k\ell}^{2}\dfrac{\partial}{\partial z}\left[\mathbf{e}_{k}^{\top}\mathbf{M}\mathbb{R}(z)\mathbf{M}^{-1}\mathbf{e}_{\ell}\cdot \mathbf{e}_{\ell}^{\top}\mathbf{M}\mathbb{R}(z)\mathbf{M}^{-1}\mathbf{e}_{k}\right]\nonumber
\end{align}
and
\begin{align}\label{eq3.4}
\textup{Cov}(M(z_{1}),M(z_{2}))=&2\left\{\dfrac{\underline{s}'(z_{1})\underline{s}'(z_{2})}{\left[\underline{s}(z_{2})-\underline{s}(z_{1})\right]^{2}}-\dfrac{1}{(z_{1}-z_{2})^{2}}\right\}\\
+&\dfrac{1}{2}\displaystyle\lim_{n \to \infty}\dfrac{1}{n} \displaystyle\sum_{k,\ell=1}^{p}r_{k\ell}^{2}\dfrac{\partial}{\partial z_{1}}\left[\mathbf{e}_{k}^{\top} \mathbf{M}\mathbb{R}(z_{1}) \mathbf{M}^{-1}\mathbf{e}_{k}\right]\nonumber\\
&\cdot \dfrac{\partial}{\partial z_{2}}\left[\mathbf{e}_{\ell}^{\top}\mathbf{M}\mathbb{R}(z_{2})\mathbf{M}^{-1}\mathbf{e}_{\ell}\right]\nonumber\\
+&\dfrac{1}{2}\displaystyle\lim_{n \to \infty} \dfrac{1}{n} \displaystyle\sum_{k,\ell=1}^{p}r_{k\ell}^{2}\dfrac{\partial}{\partial z_{1}}\left[\mathbf{e}_{k}^{\top}\mathbb{R}(z_{1})\mathbf{e}_{k}\right]\cdot \dfrac{\partial}{\partial z_{2}}\left[\mathbf{e}_{\ell}^{\top}\mathbb{R}(z_{2})\mathbf{e}_{\ell}\right]\nonumber\\
+&\dfrac{1}{2}\displaystyle\lim_{n \to \infty} \dfrac{1}{n} \displaystyle\sum_{k,\ell=1}^{p}r_{k\ell}^{2}\dfrac{\partial}{\partial z_{1}}\left[\mathbf{e}_{k}^{\top}\mathbf{M}\mathbb{R}(z_{1})\mathbf{M}^{-1}\mathbf{e}_{k}\right]\cdot \dfrac{\partial}{\partial z_{2}}\left[\mathbf{e}_{\ell}^{\top}\mathbb{R}(z_{2})\mathbf{e}_{\ell}\right]\nonumber\\
+&\dfrac{1}{2}\displaystyle\lim_{n \to \infty} \dfrac{1}{n} \displaystyle\sum_{k,\ell=1}^{p}r_{k\ell}^{2}\dfrac{\partial}{\partial z_{2}}\left[\mathbf{e}_{k}^{\top}\mathbf{M}\mathbb{R}(z_{2})\mathbf{M}^{-1}\mathbf{e}_{k}\right]\cdot \dfrac{\partial}{\partial z_{1}}\left[\mathbf{e}_{\ell}^{\top}\mathbb{R}(z_{1})\mathbf{e}_{\ell}\right]\nonumber\\
-&\underline{s}'(z_{1})\underline{s}'(z_{2})\displaystyle\lim_{n \to \infty}\dfrac{1}{n}\sum_{k=1}^{p}\left[\mathbf{e}_{k}^{\top}\mathbb{R}^{2}(z_{2})\mathbf{RMR}\mathbf{e}_{k}\right]\cdot \left[\mathbf{e}_{k}^{\top}\mathbb{R}^{2}(z_{1})\mathbf{RM}\mathbf{e}_{k}\right]\nonumber\\
-&\underline{s}'(z_{1})\underline{s}'(z_{2})\displaystyle\lim_{n \to \infty}\dfrac{1}{n}\sum_{k=1}^{p}\left[ \mathbf{e}_{k}^{\top}\mathbb{R}^{2}(z_{1})\mathbf{RMR}\mathbf{e}_{k}\right]\cdot\left[\mathbf{e}_{k}^{\top}\mathbb{R}^{2}(z_{2})\mathbf{RM}\mathbf{e}_{k}\right]\nonumber\\
-&\underline{s}'(z_{1})\underline{s}'(z_{2})\displaystyle\lim_{n \to \infty}\dfrac{1}{n}\sum_{k=1}^{p}\left[ \mathbf{e}_{k}^{\top}\mathbb{R}^{2}(z_{2})\mathbf{RMR}\mathbf{e}_{k}\right]\cdot\left[\mathbf{e}_{k}^{\top}\mathbf{M}\mathbb{R}^{2}(z_{1})\mathbf{R}\mathbf{e}_{k}\right]\nonumber\\
-&\underline{s}'(z_{1})\underline{s}'(z_{2})\displaystyle\lim_{n \to \infty}\dfrac{1}{n}\sum_{k=1}^{p}\left[ \mathbf{e}_{k}^{\top}\mathbb{R}^{2}(z_{1})\mathbf{RMR}\mathbf{e}_{k}\right]\cdot\left[\mathbf{e}_{k}^{\top}\mathbf{M}\mathbb{R}^{2}(z_{2})\mathbf{R}\mathbf{e}_{k}\right]\nonumber
\end{align}
for $\ell, \ell_{1}, \ell_{2}\in \{ 1,\ldots, K\}, \mathcal{C}, \mathcal{C}_{1} $ and $\mathcal{C}_{2}$ are three contours enclosing the support $\left[a,b\right]$ of $F^{y,H}(x)$, $\mathcal{C}_{1}, \mathcal{C}_{2}$ are non-overlapping, the contour integral $\displaystyle\oint$ is anticlockwise, $r_{k\ell}$ is the $(k, \ell)th$ element of the population correlation matrix $\mathbf{R}, \mathbf{e}_{k}$ is the kth column of $p\times p$ identity matrix $\mathbf{I}_{p}$, $\underline{s}'(z)$ is the derivative of $\underline{s}(z)$ at $z$, and
\begin{equation*}
\mathbb{R}(z)=(\mathbf{I}_{p}+\underline{s}(z)\mathbf{RM})^{-1}.
\end{equation*}
\end{theorem}
\begin{remark}
    When $\mathbf{RM}=\mathbf{I}_{p},$(\ref{eq3.3}) and (\ref{eq3.4}) can be simplified as
\begin{align*}
    \textup{E}M(z)=&y\dfrac{[\underline{s}^{\prime}(z)]^{2}}{\underline{s}(z)[1+\underline{s}(z)]^{3}}+(\tau-2)[1+z\underline{s}(z)]\dfrac{\underline{s}^{\prime}(z)}{(1+\underline{s}(z))^{2}}+\dfrac{z\underline{s}(z)+1}{2z(1+\underline{s}(z))}\\
    &+\displaystyle\lim_{n\to \infty}\dfrac{\textup{tr}\left(\mathbf{R}+\mathbf{R}^{-1}\right)}{n}\cdot \dfrac{\underline{s}^{\prime}(z)(1-\underline{s}(z))}{(1+\underline{s}(z))^{3}}-\displaystyle\lim_{n\to \infty}\dfrac{\displaystyle\sum_{k,\ell=1}^{p}\left(\mathbf{e}_{k}^{\top}\mathbf{R}\mathbf{e}_{\ell}\right)^{3}\mathbf{e}_{k}^{\top}\mathbf{R}^{-1}\mathbf{e}_{\ell}}{n}
    \\&\cdot \dfrac{\underline{s}(z)}{2z(1+\underline{s}(z))^{2}}-\displaystyle\lim_{n\to \infty}\dfrac{\displaystyle\sum_{k=1}^{p}\left(\mathbf{e}_{k}^{\top}\mathbf{R}\mathbf{e}_{k}\right)^{2}+\displaystyle\sum_{k,\ell=1}^{p}\left(\mathbf{e}_{k}^{\top}\mathbf{R}\mathbf{e}_{\ell}\right)^{3}\mathbf{e}_{k}^{\top}\mathbf{R}^{-1}\mathbf{e}_{\ell}}{n}\\
    &\qquad\qquad\qquad\qquad\cdot \dfrac{\underline{s}^{\prime}(z)}{(1+\underline{s}(z))^{3}},
\end{align*}
and
\begin{align*}
   \textup{Cov}(M(z_{1}),M(z_{2}))=&2\left\{\dfrac{\underline{s}'(z_{1})\underline{s}'(z_{2})}{\left[\underline{s}(z_{2})-\underline{s}(z_{1})\right]^{2}}-\dfrac{1}{(z_{1}-z_{2})^{2}}\right\}\\
      &+\displaystyle\lim_{n\to \infty}\dfrac{2\left[\displaystyle\sum_{k,\ell=1}^{p}\left(\mathbf{e}_{k}^{\top}\mathbf{R}\mathbf{e}_{\ell}\right)^{2}-2\textup{tr}\mathbf{R}\right]}{n}\cdot\dfrac{\underline{s}^{\prime}(z_{1})\underline{s}^{\prime}(z_{2})}{(1+\underline{s}(z_{1}))^{2}(1+\underline{s}(z_{2}))^{2}}.
\end{align*}
\end{remark}
\begin{theorem}\label{th3.2}
    Under Assumption $1$, Assumptions $3-5$ and $\mathbf{RM}=\mathbf{I}_{p},$ the random vector
    $$
    (W(g_{1}),\ldots,W(g_{K}))
    $$
    converges weakly to a multivariate Gaussian random vector $(X_{g_{1}},\ldots,X_{g_{K}})$ with
    \begin{align}
&\textup{E}X_{g_{\ell}}=\displaystyle\lim_{r\to 1^{+}}\dfrac{1}{2\pi i}\oint_{|\xi|=1}g_{\ell}(|1+\sqrt{y}\xi|^{2})\left(\dfrac{\xi}{\xi^{2}-r^{-2}}-\dfrac{1}{\xi}\right)\mathrm{d}\xi\\
&~~~~~~~~~+\dfrac{\tau-2}{2\pi i} \oint_{|\xi|=1}\dfrac{g_{\ell}(|1+\sqrt{y}\xi|^{2})}{\xi^{3}}\mathrm{d}\xi\nonumber\\
&~~~~~~~~~+\dfrac{1}{2\pi i}\oint_{|\xi|=1}g_{\ell}(|1+\sqrt{y}\xi|^{2})\dfrac{\sqrt{y}(\xi^{2}-1)}{2\xi^{3}(\xi+\sqrt{y})}\mathrm{d}\xi\nonumber\\
&~~~~~~~~~-\displaystyle\lim_{n\to \infty}\dfrac{\textup{tr}\left(\mathbf{R}+\mathbf{R}^{-1}\right)}{n}\cdot\dfrac{1}{2\pi i}\oint_{|\xi|=1}g_{\ell}(|1+\sqrt{y}\xi|^{2})\dfrac{2+\sqrt{y}\xi}{y\xi^{3}}\mathrm{d}\xi\nonumber\\
-\displaystyle\lim_{n\to \infty}&\dfrac{\displaystyle\sum_{k,\ell=1}^{p} (\mathbf{e}_{k}^{\top}\mathbf{R}\mathbf{e}_{\ell})^{3}\mathbf{e}_{k}^{\top}\mathbf{R}^{-1}\mathbf{e}_{\ell}}{n}\cdot \dfrac{1}{2\pi i}\oint_{|\xi|=1}g_{\ell}(|1+\sqrt{y}\xi|^{2}
)\dfrac{\xi^{2}-1}{2\sqrt{y}\xi^{3}(\xi+\sqrt{y})}\mathrm{d}\xi\nonumber\\
+\displaystyle\lim_{n\to \infty}&\dfrac{\displaystyle\sum_{k=1}^{p}(\mathbf{e}_{k}^{\top}\mathbf{R}\mathbf{e}_{k})^{2}+\sum_{k,\ell=1}^{p}(\mathbf{e}_{k}^{\top}\mathbf{R}\mathbf{e}_{\ell})^{3}\mathbf{e}_{k}^{\top}\mathbf{R}^{-1}\mathbf{e}_{\ell}}{n}\cdot \dfrac{1}{2\pi i} \oint_{|\xi|=1}g_{\ell}(|1+\sqrt{y}\xi|^{2})\dfrac{1+\sqrt{y}\xi}{y\xi^{3}}\mathrm{d}\xi\nonumber
     \end{align}
and
\begin{align}
\textup{Cov}(X_{g_{\ell_1}},X_{g_{\ell_2}})=&\displaystyle\lim_{r\to 1^{+}}\dfrac{-1}{2\pi^{2}}\oint_{|\xi_1|=1}\oint_{|\xi_2|=1}\dfrac{g_{\ell_{1}}(|1+\sqrt{y}\xi_{1}|^{2})g_{\ell_{2}}(|1+\sqrt{y}\xi_{2}|^{2})}{(\xi_{1}-r\xi_{2})^{2}}\mathrm{d}\xi_{2}\mathrm{d}\xi_{1}\\
&-\dfrac{1}{2\pi^{2}}\displaystyle\lim_{n\to \infty}\dfrac{\displaystyle\sum_{k,\ell=1}^{p}(\mathbf{e}_{k}^{\top}\mathbf{R}\mathbf{e}_{\ell})^{2}-2\textup{tr}\mathbf{R}}{ny}\cdot \oint_{|\xi_{1}|=1}\dfrac{g_{\ell_{1}}(|1+\sqrt{y}\xi_{1}|^{2})}{\xi_{1}^{2}}\mathrm{d}\xi_{1}\nonumber\\
&~~~\cdot\oint_{|\xi_{2}|=1}\dfrac{g_{\ell_{2}}(|1+\sqrt{y}\xi_{2}|^{2})}{\xi_{2}^{2}}\mathrm{d}\xi_{2},\nonumber
\end{align}
where $\ell,\ell_{1},\ell_{2}\in \{1,\ldots,K\},$ the contour $\displaystyle\oint$ is anticlockwise, and
$$
|1+\sqrt{y}\xi|^{2}=1+\sqrt{y}\xi+\sqrt{y}\xi^{-1}+y
$$
satisfies $|\xi|=1.$
\end{theorem}

Note that the CLT of the LSS of $\widehat{\mathbf{R}}_{n}\mathbf{M}$ in Theorem \ref{th3.1}  is the same as that of $\widehat{\mathbf{R}}_{n}$ in \cite{yin2023central} for the case $\mathbf{M}=\mathbf{I}_{p}$.

\begin{example}\label{ex3.1}
    Letting $g_{\ell}(x)=x^{\ell}$ for $\ell=1,2$ , under Assumption $1$, Assumptions $3-5$ and $\mathbf{RM}=\mathbf{I}_{p}$, we have
    \begin{itemize}
        \item Centering terms:
\begin{equation}
\begin{aligned}
      &\displaystyle\int g_{1}(x)f^{y_{n-1}}(x)\mathrm{d}x=1,\\
      &\displaystyle\int g_{2}(x)f^{y_{n-1}}(x)\mathrm{d}x=1+y_{n-1}.
\end{aligned}
\end{equation}
      \item Mean terms:
     \begin{equation}
\begin{aligned}
      \textup{E}X_{g_{1}}=&\dfrac{3y}{2}-a_{\mathbf{R}}+\dfrac{b_{\mathbf{R}}}{2},\\
      \textup{E}X_{g_{2}}=&\dfrac{3}{2}y^2+(\tau+5-2a_{\mathbf{R}}+\dfrac{5}{2}b_{\mathbf{R}})y-(1+4a_{\mathbf{R}})+b_{\mathbf{R}}y^{-1},
\end{aligned}
\end{equation}
\item Variance and covariance terms:
\begin{equation}
    \begin{aligned}
        &\textup{Var}(X_{g_1})=2c_{\mathbf{R}}-2y,\\
        &\textup{Var}(X_{g_2})=4y^2+8(1+y)^{2}(c_{\mathbf{R}}-y),\\
       &\textup{Cov}(X_{g_1},X_{g_2})=4(1+y)(c_{\mathbf{R}}-y),
    \end{aligned}
\end{equation}
    \end{itemize}
where $a_{\mathbf{R}}=\displaystyle\lim_{n\to \infty} \dfrac{\textup{tr}(\mathbf{R}+\mathbf{R}^{-1})}{n}, b_{\mathbf{R}}=\displaystyle\lim_{n\to \infty} \dfrac{\displaystyle\sum_{k,\ell=1}^{p}(\mathbf{e}_{k}^{\top}\mathbf{R}\mathbf{e}_{\ell})^{3}\mathbf{e}_{k}^{\top}\mathbf{R}^{-1}\mathbf{e}_{\ell}}{n}, c_{\mathbf{R}}=\displaystyle\lim_{n\to \infty}\\\dfrac{\displaystyle\sum_{k,\ell=1}^{p}(\mathbf{e}_{k}^{\top}\mathbf{R}\mathbf{e}_{\ell})^{2}}{n}.$
\end{example}

\section{Central limit theorem of linear spectral statistics of $\widehat{\mathbf{R}}_{n}\mathbf{M}$ under linear independent component structure}
This section establishes the CLT of the random vector~$(W(g_{1}),\ldots,W(g_{K}))$~under the linear independent component structure assumption $2$.

\begin{theorem}\label{th4.1}
Under Assumptions $A_{L}$-B-C-D,
the random vector $(W(g_{1}),\ldots\\
, W(g_{K}))$ weakly converges to a multivariate Gaussian vector $(X_{g_{1}},\ldots, X_{g_{K}})$ with the mean and covariance functions as follows:
\begin{equation}\label{temp3}
\textup{E}X_{g_{\ell}}=-\dfrac{1}{2\pi i}\displaystyle\oint_{\mathcal{C}}g_{\ell}(z)\textup{E}M(z)\textup{d}z,
\end{equation}
\begin{equation}\label{temp4}
\textup{Cov}(X_{g_{\ell_{1}}},X_{g_{\ell_{2}}})=-\dfrac{1}{4\pi^{2}}\displaystyle\oint_{\mathcal{C}_{1}}\oint_{\mathcal{C}_{2}}g_{\ell_{1}}(z_{1})g_{\ell_{2}}(z_{2})\textup{Cov}(M(z_{1}),M(z_{2}))\textup{d}z_{1}\textup{d}z_{2},
\end{equation}
where
\begin{align}\label{eq4.1}
\textup{E}M(z)=&y\displaystyle\int\dfrac{\left[t\underline{s}'(z)\right]^{2}}{\underline{s}(z)\left[1+t\underline{s}(z)\right]^{3}}\mathrm{d}H(t)\\
&+\beta_x y \underline{s}(z)\underline{s}'(z)\lim_{n\to\infty}\dfrac{1}{n}\sum_{k=1}^{p}\mathbf{e}_{k}^\top\mathbf{G}^\top\mathbf{M}\mathbb{R}(z)\mathbf{G}\mathbf{e}_{k}\cdot\mathbf{e}_{k}^\top\mathbf{G}^\top\mathbf{M}\mathbb{R}^2(z)\mathbf{G}\mathbf{e}_{k}\nonumber\\
&+\displaystyle\lim_{n\to \infty } \dfrac{1}{n} \sum_{k=1}^{p}\dfrac{\partial}{\partial z}\left[\underline{s}(z)\mathbf{e}_{k}^{\top}\mathbf{RM}\mathbb{R}(z)\mathbf{e}_{k}\cdot \mathbf{e}_{k}^{\top}\mathbb{R}(z)\mathbf{R}\mathbf{e}_{k}\right]\nonumber\\
&+\displaystyle\lim_{n\to\infty} \dfrac{1}{n} \sum_{k=1}^{p}\dfrac{\partial}{\partial z}\left[\underline{s}(z)\mathbf{e}_{k}^{\top}\mathbf{RM}\mathbb{R}(z)\mathbf{M}\mathbf{e}_{k}\cdot \mathbf{e}_{k}^{\top}\mathbf{M}\mathbb{R}(z)\mathbf{R}\mathbf{e}_{k}\right]\nonumber\\
&+\displaystyle\lim_{n\to \infty } \dfrac{\beta_{x}}{2n} \sum_{k,\ell=1}^{p}g_{k\ell}^{2}\cdot\dfrac{\partial}{\partial z}\left[\underline{s}(z)\mathbf{e}_{\ell}^{\top}\mathbf{G}^{\top}\mathbf{M}\mathbb{R}(z)\mathbf{e}_{k}\cdot\mathbf{e}_{k}^{\top}\mathbb{R}(z)\mathbf{Ge}_{\ell}\right]\nonumber\\
&+\displaystyle\lim_{n\to \infty } \dfrac{\beta_{x}}{2n} \sum_{k,\ell=1}^{p}g_{k\ell}^{2}\cdot\dfrac{\partial}{\partial z}\left[\underline{s}(z)\mathbf{e}_{\ell}^{\top}\mathbf{G}^{\top}\mathbf{M}\mathbb{R}(z)\mathbf{M}\mathbf{e}_{k}\cdot\mathbf{e}_{k}^{\top}\mathbf{M}\mathbb{R}(z)\mathbf{Ge}_{\ell}\right]\nonumber\\
&+\dfrac{1}{4z}\displaystyle\lim_{n\to \infty}\dfrac{1}{n}\sum_{\ell=1}^{p}\left[\mathbf{e}_{\ell}^{\top}\mathbb{R}(z)\mathbf{e}_{\ell}\right]\cdot\left[\dfrac{\beta_{x}}{2}\sum_{j=1}^{p}g_{\ell j}^{4}+1\right]\nonumber\\
&+\dfrac{1}{4z}\displaystyle\lim_{n\to \infty}\dfrac{1}{n}\sum_{\ell=1}^{p}\left[\mathbf{e}_{\ell}^{\top}\mathbf{M}\mathbb{R}(z)\mathbf{M}^{-1}\mathbf{e}_{\ell}\right]\cdot\left[\dfrac{\beta_{x}}{2}\sum_{j=1}^{p}g_{\ell j}^{4}+1\right]\nonumber\\
&-\dfrac{1}{4z}\displaystyle\lim_{n\to \infty}\dfrac{1}{n}\sum_{\ell=1}^{p}\left[\mathbf{e}_{\ell}^{\top}\mathbb{R}^{2}(z)\mathbf{e}_{\ell}\right]\cdot\left[\dfrac{\beta_{x}}{2}\sum_{j=1}^{p}g_{\ell j}^{4}+1\right]\nonumber\\
&-\dfrac{1}{4z}\displaystyle\lim_{n\to \infty}\dfrac{1}{n}\sum_{\ell=1}^{p}\left[\mathbf{e}_{\ell}^{\top}\mathbf{M}\mathbb{R}^{2}(z)\mathbf{M}^{-1}\mathbf{e}_{\ell}\right]\cdot\left[\dfrac{\beta_{x}}{2}\sum_{j=1}^{p}g_{\ell j}^{4}+1\right]\nonumber\\
&-\dfrac{1}{2z}\displaystyle\lim_{n\to \infty}\dfrac{1}{n}\sum_{k,\ell=1}^{p}\mathbf{e}_{\ell}^{\top}\mathbf{M}^{-1}\mathbf{e}_{k}\left[\mathbf{e}_{k}^{\top}\mathbf{M}\mathbb{R}(z)\mathbf{e}_{\ell}\right]\cdot\left[r_{k\ell}^{2}+\dfrac{\beta_{x}}{2}\sum_{j=1}^{p}g^{2}_{\ell j}g^{2}_{kj}\right]\nonumber\\
&+\dfrac{1}{2z}\displaystyle\lim_{n\to \infty}\dfrac{1}{n}\sum_{k,\ell=1}^{p}\mathbf{e}_{\ell}^{\top}\mathbf{M}^{-1}\mathbf{e}_{k}\left[\mathbf{e}_{k}^{\top}\mathbf{M}\mathbb{R}^{2}(z)\mathbf{e}_{\ell}\right]\cdot\left[r_{k\ell}^{2}+\dfrac{\beta_{x}}{2}\sum_{j=1}^{p}g^{2}_{\ell j}g^{2}_{kj}\right]\nonumber\\
&+\dfrac{1}{2}\displaystyle\lim_{n\to \infty}\dfrac{1}{n}\sum_{k,\ell=1}^{p}\dfrac{\partial}{\partial z}\left[\mathbf{e}_{k}^{\top}\mathbb{R}(z)\mathbf{M}^{-1}\mathbf{e}_{\ell}\cdot \mathbf{e}_{\ell}^{\top}\mathbf{M}\mathbb{R}(z)\mathbf{e}_{k}\right]\cdot\left[r_{k\ell}^{2}+\dfrac{\beta_{x}}{2}\sum_{j=1}^{p}g^{2}_{\ell j}g^{2}_{kj}\right]\nonumber\\
&+\dfrac{1}{4}\displaystyle\lim_{n\to \infty}\dfrac{1}{n}\sum_{k,\ell=1}^{p}\dfrac{\partial}{\partial z}\left[\mathbf{e}_{k}^{\top}\mathbb{R}(z)\mathbf{e}_{\ell}\cdot \mathbf{e}_{\ell}^{\top}\mathbb{R}(z)\mathbf{e}_{k}\right]\cdot\left[r_{k\ell}^{2}+\dfrac{\beta_{x}}{2}\sum_{j=1}^{p}g^{2}_{\ell j}g^{2}_{kj}\right]\nonumber\\
&+\dfrac{1}{4}\displaystyle\lim_{n\to \infty}\dfrac{1}{n}\sum_{k,\ell=1}^{p}\dfrac{\partial}{\partial z}\left[\mathbf{e}_{k}^{\top}\mathbf{M}\mathbb{R}(z)\mathbf{M}^{-1}\mathbf{e}_{\ell}\cdot \mathbf{e}_{\ell}^{\top}\mathbf{M}\mathbb{R}(z)\mathbf{M}^{-1}\mathbf{e}_{k}\right]\nonumber\\
&~~\cdot\left[r_{k\ell}^{2}+\dfrac{\beta_{x}}{2}\sum_{j=1}^{p}g^{2}_{\ell j}g^{2}_{kj}\right]\nonumber
\end{align}
and
\begin{align}
&\textup{Cov}(M(z_{1}),M(z_{2}))\\
=&2\left\{\dfrac{\underline{s}'(z_{1})\underline{s}'(z_{2})}{\left[\underline{s}(z_{2})-\underline{s}(z_{1})\right]^{2}}-\dfrac{1}{(z_{1}-z_{2})^{2}}\right\}\nonumber\\
&+\lim_{n\to \infty}\dfrac{\beta_x y}{n}\sum_{k=1}^p\dfrac{\partial}{\partial z_{1}}\underline{s}(z_{1})\mathbf{e}_{k}^{\top}\mathbf{G}^\top\mathbf{M}\mathbb{R}(z_1)\mathbf{G}\mathbf{e}_k\cdot\dfrac{\partial}{\partial z_2}\underline{s}(z_2)\mathbf{e}_{k}^{\top}\mathbf{G}^\top\mathbf{M}\mathbb{R}(z_2)\mathbf{G}\mathbf{e}_k\nonumber\\
&+\dfrac{1}{4}\displaystyle\lim_{n\to \infty}\dfrac{1}{n}\sum_{k,\ell=1}^{p}\left(\beta_{x}\sum_{i=1}^{p}g_{ki}^{2}g_{\ell i}^{2}+2r_{k\ell}^{2}\right)\cdot\dfrac{\partial}{\partial z_{1}}\mathbf{e}_{k}^{\top}\mathbb{R}(z_{1})\mathbf{e}_{k}\cdot\dfrac{\partial}{\partial z_{2}}\mathbf{e}_{\ell}^{\top}\mathbb{R}(z_{2})\mathbf{e}_{\ell}\nonumber\\
&+\dfrac{1}{4}\displaystyle\lim_{n\to \infty}\dfrac{1}{n}\sum_{k,\ell=1}^{p}\left(\beta_{x}\sum_{i=1}^{p}g_{ki}^{2}g_{\ell i}^{2}+2r_{k\ell}^{2}\right)\cdot\dfrac{\partial}{\partial z_{1}}\mathbf{e}_{k}^{\top}\mathbf{M}\mathbb{R}(z_{1})\mathbf{M}^{-1}\mathbf{e}_{k}\nonumber\\
&~~\cdot\dfrac{\partial}{\partial z_{2}}\mathbf{e}_{\ell}^{\top}\mathbf{M}\mathbb{R}(z_{2})\mathbf{M}^{-1}\mathbf{e}_{\ell}\nonumber\\
&+\dfrac{1}{4}\displaystyle\lim_{n\to \infty}\dfrac{1}{n}\sum_{k,\ell=1}^{p}\left(\beta_{x}\sum_{i=1}^{p}g_{ki}^{2}g_{\ell i}^{2}+2r_{k\ell}^{2}\right)\dfrac{\partial}{\partial z_{1}}\mathbf{e}_{k}^{\top}\mathbb{R}(z_{1})\mathbf{e}_{k}\cdot\dfrac{\partial}{\partial z_{2}}\mathbf{e}_{\ell}^{\top}\mathbf{M}\mathbb{R}(z_{2})\mathbf{M}^{-1}\mathbf{e}_{\ell}\nonumber\\
&+\dfrac{1}{4}\displaystyle\lim_{n\to \infty}\dfrac{1}{n}\sum_{k,\ell=1}^{p}\left(\beta_{x}\sum_{i=1}^{p}g_{ki}^{2}g_{\ell i}^{2}+2r_{k\ell}^{2}\right)\dfrac{\partial}{\partial z_{2}}\mathbf{e}_{k}^{\top}\mathbb{R}(z_{2})\mathbf{e}_{k}\cdot\dfrac{\partial}{\partial z_{1}}\mathbf{e}_{\ell}^{\top}\mathbf{M}\mathbb{R}(z_{1})\mathbf{M}^{-1}\mathbf{e}_{\ell}\nonumber\\
&-\displaystyle\lim_{n\to \infty}\dfrac{1}{n}\sum_{k=1}^{p}\dfrac{\partial}{\partial z_{1}}\mathbf{e}_{k}^{\top}\mathbb{R}(z_{1})\mathbf{e}_{k}\cdot \dfrac{\partial}{\partial z_{2}}\mathbf{e}_{k}^{\top}\mathbb{R}(z_{2})\mathbf{R}\mathbf{e}_{k}\nonumber\\
&-\displaystyle\lim_{n\to \infty}\dfrac{1}{n}\sum_{k=1}^{p}\dfrac{\partial}{\partial z_{2}}\mathbf{e}_{k}^{\top}\mathbb{R}(z_{2})\mathbf{e}_{k}\cdot \dfrac{\partial}{\partial z_{1}}\mathbf{e}_{k}^{\top}\mathbb{R}(z_{1})\mathbf{R}\mathbf{e}_{k}\nonumber\\
&+\displaystyle\lim_{n\to \infty}\dfrac{\beta_{x}}{2n}\sum_{k,\ell=1}^{p}g_{k\ell}^{2}\dfrac{\partial}{\partial z_{1}}\mathbf{e}_{k}^{\top}\mathbb{R}(z_{1})\mathbf{e}_{k}\cdot \dfrac{\partial}{\partial z_{2}}\underline{s}(z_{2})\mathbf{e}_{\ell}^{\top}\mathbf{G}^{\top}\mathbb{R}(z_{2})\mathbf{G}\mathbf{e}_{\ell}\nonumber\\
&+\displaystyle\lim_{n\to \infty}\dfrac{\beta_{x}}{2n}\sum_{k,\ell=1}^{p}g_{k\ell}^{2}\dfrac{\partial}{\partial z_{2}}\mathbf{e}_{k}^{\top}\mathbb{R}(z_{2})\mathbf{e}_{k}\cdot \dfrac{\partial}{\partial z_{1}}\underline{s}(z_{1})\mathbf{e}_{\ell}^{\top}\mathbf{G}^{\top}\mathbb{R}(z_{1})\mathbf{G}\mathbf{e}_{\ell}\nonumber\\
&-\displaystyle\lim_{n\to \infty}\dfrac{1}{n}\sum_{k=1}^{p}\dfrac{\partial}{\partial z_{1}}\mathbf{e}_{k}^{\top}\mathbf{M}\mathbb{R}(z_{1})\mathbf{M}^{-1}\mathbf{e}_{k}\cdot \dfrac{\partial}{\partial z_{2}}\mathbf{e}_{k}^{\top}\mathbb{R}(z_{2})\mathbf{R}\mathbf{e}_{k}\nonumber\\
&-\displaystyle\lim_{n\to \infty}\dfrac{1}{n}\sum_{k=1}^{p}\dfrac{\partial}{\partial z_{2}}\mathbf{e}_{k}^{\top}\mathbf{M}\mathbb{R}(z_{2})\mathbf{M}^{-1}\mathbf{e}_{k}\cdot \dfrac{\partial}{\partial z_{1}}\mathbf{e}_{k}^{\top}\mathbb{R}(z_{1})\mathbf{R}\mathbf{e}_{k}\nonumber\\
&+\displaystyle\lim_{n\to \infty}\dfrac{\beta_{x}}{2n}\sum_{k,\ell=1}^{p}g_{k\ell}^{2}\dfrac{\partial}{\partial z_{1}}\mathbf{e}_{k}^{\top}\mathbf{M}\mathbb{R}(z_{1})\mathbf{M}^{-1}\mathbf{e}_{k}\cdot \dfrac{\partial}{\partial z_{2}}\underline{s}(z_{2})\mathbf{e}_{\ell}^{\top}\mathbf{G}^{\top}\mathbb{R}(z_{2})\mathbf{G}\mathbf{e}_{\ell}\nonumber\\
&+\displaystyle\lim_{n\to \infty}\dfrac{\beta_{x}}{2n}\sum_{k,\ell=1}^{p}g_{k\ell}^{2}\dfrac{\partial}{\partial z_{2}}\mathbf{e}_{k}^{\top}\mathbf{M}\mathbb{R}(z_{2})\mathbf{M}^{-1}\mathbf{e}_{k}\cdot \dfrac{\partial}{\partial z_{1}}\underline{s}(z_{1})\mathbf{e}_{\ell}^{\top}\mathbf{G}^{\top}\mathbb{R}(z_{1})\mathbf{G}\mathbf{e}_{\ell},\nonumber
\end{align}
for $\ell, \ell_{1}, \ell_{2}\in \{ 1,\ldots, K\}, \mathcal{C}, \mathcal{C}_{1} $ and $\mathcal{C}_{2}$ are three contours enclosing the support $\left[a,b\right]$ of $F^{y,H}(x)$, $\mathcal{C}_{1}, \mathcal{C}_{2}$ are non-overlapping, the contour integral $\displaystyle\oint$ is anticlockwise, $r_{k\ell}$ is the $(k, l)th$ element of the population correlation matrix $\mathbf{R}$,~$g_{k\ell}$ is the $(k, l)th$ element of ~$\mathbf{G}$,~$\mathbf{e}_{k}$ is the kth column of $p\times p$ identity matrix $\mathbf{I}_{p}$, $\underline{s}'(z)$ is the derivative of $\underline{s}(z)$ at $z$,~and
\begin{equation*}
\mathbb{R}(z)=(\mathbf{I}_{p}+\underline{s}(z)\mathbf{RM})^{-1}.
\end{equation*}
\end{theorem}

Note that the CLT of the LSS of $\widehat{\mathbf{R}}_{n}\mathbf{M}$ in Theorem \ref{th4.1}  is as same as that of $\widehat{\mathbf{R}}_{n}$ in \cite{yin2023central} for the case $\mathbf{M}=\mathbf{I}_{p}$.

\section{Joint limiting distribution of linear spectral statistics of $\widehat{\mathbf{R}}_{n}\mathbf{M}_{1},\ldots\\,\widehat{\mathbf{R}}_{n}\mathbf{M}_{K}$}

Let $\mathbf{M}_1,\ldots,\mathbf{M}_K$ be the positive-definite matrices, the ESD of $\mathbf{RM}_j$ as $H_{nj}$ and the corresponding LSD is $H_j$. Now define the LSS as
\begin{equation}
\begin{aligned}
     L_{g_{\ell j}}&=\displaystyle\sum_{i=1}^{p}g_{\ell j}(\widehat{\lambda}^j_{i}),~~~\ell=1,\ldots,m_j, j=1,\ldots,K
\end{aligned}
\end{equation}
where $g_{\ell j}(\cdot)$ are some known analytic functions, $\widehat{\lambda}^j_{i}$ is the $i$th largest eigenvalue of $\widehat{\mathbf{R}}_{n}\mathbf{M}_{j}$.
 Let the $m_1+\ldots+m_K$ dimensional random vector be
 $(W(g_{11}),\ldots,W(g_{m_11}),\\\ldots
 ,W(g_{1K}),\ldots,W(g_{m_KK})$
,  where
\begin{equation}
   W(g_{\ell j})=\sum_{i=1}^{p}g_{\ell j}(\widehat{\lambda}^j_{i})\delta(\{\widehat{\lambda}^j_{i}>0\})-p\displaystyle\int_{a}^{b}g_{\ell j}(x)f^{y_{n-1},H_{nj}}(x)\mathrm{d}x,
\end{equation}
and $f^{y_{n-1},H_{nj}}$ is defined in (\ref{eq2.6}) with $y_{n-1}=p/(n-1)$ and $H_{n}$ being replaced by $H_{nj}$.

\subsection{Joint limiting distribution of LSSs of $\widehat{\mathbf{R}}_{n}\mathbf{M}_{1}$, $\ldots$, $\widehat{\mathbf{R}}_{n}\mathbf{M}_{K}$ under elliptical structure}

The following theorem will give the central limit theorem of the random vector $(W(g_{11})$, $\ldots$, $W(g_{m_11})$, $\ldots$, $W(g_{1K})$, $\ldots$, $W(g_{m_KK})$.

\begin{theorem}
    Under Assumptions $1, 3, 4, 5$, the random vector $(W(g_{11})$, $\ldots$, $W(g_{m_11})$, $\ldots$, $W(g_{1K})$, $\ldots$, $W(g_{m_KK})$ weakly converges to a multivariate Gaussian vector $(X_{g_{11}}$, $\ldots$, $X_{g_{m_11}}$, $\ldots$, $X_{g_{1K}}$, $\ldots$, $X_{g_{m_KK}}$ with mean function and covariance functions as follows:
\begin{itemize}
\item
(1). ${\rm E}X_{g_{\ell j}}$ is in (\ref{temp3}) by replacing $g_{\ell}$ by $g_{\ell j}$;
\item
(2). ${\rm Cov}(X_{g_{\ell_1 j}}, X_{g_{\ell_2 j}})$ is in (\ref{temp4}) by replacing $g_{\ell_1}$ and $g_{\ell_2}$ by $g_{\ell_1 j}$ and $g_{\ell_2 j}$;
\item
(3).
${\rm Cov}(X_{g_{\ell_1 j}}, X_{g_{\ell_2 h}})=-\dfrac{1}{4\pi^{2}}\displaystyle\oint_{\mathcal{C}_{1}}\oint_{\mathcal{C}_{2}}g_{\ell_1 j}(z_{1})g_{\ell_2 h}(z_{2})\textup{Cov}(M^j(z_{1}),M^h(z_{2}))$\\
$\mathrm{d}z_{2}\mathrm{d}z_{1}$, where

\begin{align}\label{eq5.3}
&\textup{Cov}(M^j(z_{1}),M^h(z_{2}))\\
=&2\dfrac{\partial^2}{\partial z_1 \partial z_2}\log(1- a(z_{1},z_{2}))+\lim_{n\to \infty}\dfrac{\beta_x y}{n}\sum_{k=1}^p\dfrac{\partial}{\partial z_{1}}\underline{s}_j(z_1)\mathbf{e}_{k}^{\top}\mathbf{G}^\top\mathbf{M}_j\mathbb{R}_j(z_1)\mathbf{G}\mathbf{e}_k\nonumber\\
&\cdot\dfrac{\partial}{\partial z_2}\underline{s}_h(z_2)\mathbf{e}_{k}^{\top}\mathbf{G}^\top\mathbf{M}_h\mathbb{R}_h(z_2)\mathbf{G}\mathbf{e}_k\nonumber\\
+&\dfrac{1}{4}\displaystyle\lim_{n\to \infty}n^{-1}\sum_{k,\ell=1}^{p}\left(\beta_{x}\sum_{i=1}^{p}g_{ki}^{2}g_{\ell i}^{2}+2r_{k\ell}^{2}\right)\cdot\dfrac{\partial}{\partial z_{1}}\mathbf{e}_{k}^{\top}\mathbb{R}_j(z_1)\mathbf{e}_{k}\cdot\dfrac{\partial}{\partial z_{2}}\mathbf{e}_{\ell}^{\top}\mathbb{R}_h(z_2)\mathbf{e}_{\ell}\nonumber\\
+&\dfrac{1}{4}\displaystyle\lim_{n\to \infty}n^{-1}\sum_{k,\ell=1}^{p}\left(\beta_{x}\sum_{i=1}^{p}g_{ki}^{2}g_{\ell i}^{2}+2r_{k\ell}^{2}\right)\cdot\dfrac{\partial}{\partial z_{1}}\mathbf{e}_{k}^{\top}\mathbf{M}_j\mathbb{R}_j(z_1)\mathbf{M}_j^{-1}\mathbf{e}_{k}\nonumber\\
&\cdot\dfrac{\partial}{\partial z_{2}}\mathbf{e}_{\ell}^{\top}\mathbf{M}_h\mathbb{R}_h(z_2)\mathbf{M}_h^{-1}\mathbf{e}_{\ell}\nonumber\\
+&\dfrac{1}{4}\displaystyle\lim_{n\to \infty}n^{-1}\sum_{k,\ell=1}^{p}\left(\beta_{x}\sum_{i=1}^{p}g_{ki}^{2}g_{\ell i}^{2}+2r_{k\ell}^{2}\right)\dfrac{\partial}{\partial z_{1}}\mathbf{e}_{k}^{\top}\mathbb{R}_j(z_1)\mathbf{e}_{k}\nonumber\\
&\cdot\dfrac{\partial}{\partial z_{2}}\mathbf{e}_{\ell}^{\top}\mathbf{M}_h\mathbb{R}_h(z_2)\mathbf{M}_h^{-1}\mathbf{e}_{\ell}\nonumber\\
+&\dfrac{1}{4}\displaystyle\lim_{n\to \infty}n^{-1}\sum_{k,\ell=1}^{p}\left(\beta_{x}\sum_{i=1}^{p}g_{ki}^{2}g_{\ell i}^{2}+2r_{k\ell}^{2}\right)\dfrac{\partial}{\partial z_{1}}\mathbf{e}_{\ell}^{\top}\mathbf{M}_j\mathbb{R}_j(z_1)\mathbf{M}_j^{-1}\mathbf{e}_{\ell}\nonumber\\
&\cdot\dfrac{\partial}{\partial z_{2}}\mathbf{e}_{k}^{\top}\mathbb{R}_h(z_2)\mathbf{e}_{k}\nonumber\\
-&\displaystyle\lim_{n\to \infty}n^{-1}\sum_{k=1}^{p}\dfrac{\partial}{\partial z_{1}}\mathbf{e}_{k}^{\top}\mathbb{R}_j(z_1)\mathbf{e}_{k}\cdot \dfrac{\partial}{\partial z_{2}}\mathbf{e}_{k}^{\top}\mathbb{R}_h(z_2)\mathbf{R}\mathbf{e}_{k}\nonumber\\
-&\displaystyle\lim_{n\to \infty}n^{-1}\sum_{k=1}^{p}\dfrac{\partial}{\partial z_{1}}\mathbf{e}_{k}^{\top}\mathbb{R}_j(z_1)\mathbf{R}\mathbf{e}_{k}\cdot\dfrac{\partial}{\partial z_{2}}\mathbf{e}_{k}^{\top}\mathbb{R}_h(z_2)\mathbf{e}_{k}\nonumber\\
+&\displaystyle\lim_{n\to \infty}\dfrac{\beta_{x}}{2n}\sum_{k,\ell=1}^{p}g_{k\ell}^{2}\dfrac{\partial}{\partial z_{1}}\mathbf{e}_{k}^{\top}\mathbb{R}_j(z_1)\mathbf{e}_{k}\cdot \dfrac{\partial}{\partial z_{2}}\underline{s}(z_{2})\mathbf{e}_{\ell}^{\top}\mathbf{G}^{\top}\mathbb{R}_h(z_2)\mathbf{G}\mathbf{e}_{\ell}\nonumber\\
+&\displaystyle\lim_{n\to \infty}\dfrac{\beta_{x}}{2n}\sum_{k,\ell=1}^{p}g_{k\ell}^{2}\dfrac{\partial}{\partial z_{1}}\underline{s}(z_{1})\mathbf{e}_{\ell}^{\top}\mathbf{G}^{\top}\mathbb{R}_j(z_1)\mathbf{G}\mathbf{e}_{\ell}\cdot\dfrac{\partial}{\partial z_{2}}\mathbf{e}_{k}^{\top}\mathbb{R}_h(z_2)\mathbf{e}_{k}\nonumber\\
-&\displaystyle\lim_{n\to \infty}n^{-1}\sum_{k=1}^{p}\dfrac{\partial}{\partial z_{1}}\mathbf{e}_{k}^{\top}\mathbf{M}_j\mathbb{R}_j(z_1)\mathbf{M}_j^{-1}\mathbf{e}_{k}\cdot \dfrac{\partial}{\partial z_{2}}\mathbf{e}_{k}^{\top}\mathbb{R}_h(z_2)\mathbf{R}\mathbf{e}_{k}\nonumber\\
-&\displaystyle\lim_{n\to \infty}n^{-1}\sum_{k=1}^{p}\dfrac{\partial}{\partial z_{1}}\mathbf{e}_{k}^{\top}\mathbb{R}_j(z_1)\mathbf{R}\mathbf{e}_{k}\cdot\dfrac{\partial}{\partial z_{2}}\mathbf{e}_{k}^{\top}\mathbf{M}_h\mathbb{R}_h(z_2)\mathbf{M}_h^{-1}\mathbf{e}_{k}\nonumber\\
+&\displaystyle\lim_{n\to \infty}\dfrac{\beta_{x}}{2n}\sum_{k,\ell=1}^{p}g_{k\ell}^{2}\dfrac{\partial}{\partial z_{1}}\mathbf{e}_{k}^{\top}\mathbf{M}_j\mathbb{R}_j(z_1)\mathbf{M}_j^{-1}\mathbf{e}_{k}\cdot \dfrac{\partial}{\partial z_{2}}\underline{s}_h(z_{2})\mathbf{e}_{\ell}^{\top}\mathbf{G}^{\top}\mathbb{R}_h(z_2)\mathbf{G}\mathbf{e}_{\ell}\nonumber\\
+&\displaystyle\lim_{n\to \infty}\dfrac{\beta_{x}}{2n}\sum_{k,\ell=1}^{p}g_{k\ell}^{2}\dfrac{\partial}{\partial z_{1}}\underline{s}_j(z_{1})\mathbf{e}_{\ell}^{\top}\mathbf{G}^{\top}\mathbb{R}_j(z_1)\mathbf{G}\mathbf{e}_{\ell}\cdot\dfrac{\partial}{\partial z_{2}}\mathbf{e}_{k}^{\top}\mathbf{M}_h\mathbb{R}_h(z_2)\mathbf{M}_h^{-1}\mathbf{e}_{k},\nonumber
\end{align}
\end{itemize}
for $ \ell \in \{ 1,\ldots, K_1\},  j \in \{1,\ldots, K_2\}, \mathcal{C}_{1} $ and $\mathcal{C}_{2}$ are two contours enclosing the supports of $F^{y,H_j}(x)$ and $F^{y,H_h}(x)$, $\mathcal{C}_{1}, \mathcal{C}_{2}$ are non-overlapping, the contour integral $\displaystyle\oint$ is anticlockwise, $r_{k\ell}$ is the $(k, \ell)th$ element of the population correlation matrix $\mathbf{R}, \mathbf{e}_{k}$ is the kth column of $p\times p$ identity matrix $\mathbf{I}_{p}$, $\underline{s}_j'(z)$ is the derivative of $\underline{s}_j(z)$ at $z$, ${\underline{s}_h}^{\prime}(z)$ is the derivative of $\underline{s}_h(z)$ at z, $\underline{s}_j(z)$ is the Stieltjes transform of the LSD $F^{y,H_j}$ defined in (\ref{eq2.5}),
and
\begin{equation*}
\mathbb{R}_j(z_1)=[\mathbf{I}_{p}+\underline{s}_j(z_1)\mathbf{RM}_j]^{-1},\quad\mathbb{R}_h(z_2)=[\mathbf{I}_{p}+\underline{s}_h(z_2)\mathbf{RM}_h]^{-1},
\end{equation*}
\begin{equation*}
    a(z_1,z_2)=\underline{s}_j(z_{1})\underline{s}_h(z_2)\cdot \lim_{n\to \infty}n^{-1}{\rm tr}\left[\mathbf{RM}_j\mathbb{R}_j(z_1)\mathbf{RM}_h\mathbb{R}_h(z_2)\right].
\end{equation*}
\end{theorem}
\section{An application}
Testing method: We study the hypothesis testing problem that the population correlation matrix is equal to a given matrix as follows
\begin{equation}\label{eq6.1}
    H_{0}: \mathbf{R}=\mathbf{R}_{0}~~~v.s.~~~ H_{1}: \mathbf{R}\neq \mathbf{R}_{0},
\end{equation}
where $\mathbf{R}_{0}$ is a pre-specified matrix. Let $\widehat{\mathbf{R}}_{n}$ be the sample correlation matrix. Based on the difference and ratio between $\widehat{\mathbf{R}}_{n}$ and $\mathbf{R}_0$, we propose the following test statistic
$$
\max\left\{\sigma_1^{-1}|T_1-\mu_1|,~~\sigma_2^{-1}|T_2-\mu_2|\right\}
$$
where $$T_1=\textup{tr}\left(\widehat{\mathbf{R}}_{n}\mathbf{R}_0^{-1}-\mathbf{I}_p\right)^{2},~~T_2=\textup{tr}\left(\widehat{\mathbf{R}}_n-\mathbf{R}_0\right)^2.$$
$\mu_1, \mu_2, \sigma_1, \sigma_2$ will be shown in the following theorem.

\begin{theorem}
    Under Assumptions $1, 3, 4, 5$ and under $H_0$, we have
\begin{itemize}
    \item  (1). $\sigma_1^{-1}(T_1-\mu_1)\to N(0,1)$, where
   \begin{align*}
         \mu_1&=py_{n-1}-2\textup{E}X_{g_1}+\textup{E}X_{g_2},\\
         \sigma_1&=\textup{Var}(X_{g_2})+4\textup{Var}X_{g_1}-4\textup{Cov}(X_{g_1},X_{g_2}),
   \end{align*}
  where  $y_{n-1}=p/(n-1), \textup{E}X_{g_1}, \textup{E}X_{g_2}, \textup{Var}(X_{g_1}), \textup{Var}(X_{g_1}), \textup{Cov}(X_{g_1},X_{g_2})$ are in Example \ref{ex3.1} and $y$ can also be replaced by $y_{n-1}$.

  \item (2). $\sigma_2^{-1}(T_2-\mu_2)\to N(0,1)$, where
   \begin{align*}
         \mu_2&=(p+\tau-3)y_{n-1}-p+\dfrac{1}{\pi i}\displaystyle\oint_{\mathcal{C}}z\textup{E}M(z)\mathrm{d}z+\textup{tr}(\mathbf{R}_0^2),\\
         \sigma_2^2&=4y_{n-1}^2-\dfrac{1}{\pi^2}\displaystyle\oint_{\mathcal{C}_1}\oint_{\mathcal{C}_2}z_1z_2 \textup{Cov}(M^2(z_1),M^2(z_2))\mathrm{d}z_2\mathrm{d}z_1 \\
         &~~~+\dfrac{1}{\pi^{2}}\displaystyle\oint_{\mathcal{C}_1}\oint_{\mathcal{C}_2}z_1^2 z_2\textup{Cov}(M^{1}(z_1),M^{2}(z_2))\mathrm{d}z_2\mathrm{d}z_1,
   \end{align*}
   where $\textup{E}M(z)$ is defined in (\ref{eq3.3}) with $\mathbf{M}=\mathbf{R}_0$, $\textup{Cov}(M^2(z_1),M^2(z_2))$ is defined in (\ref{eq3.4}) with $\mathbf{M}=\mathbf{R}_0$ and $\textup{Cov}(M^{1}(z_1),M^{2}(z_2))$ is defined in (\ref{eq5.3}) with $\mathbf{M}_1=\mathbf{I}_p, \mathbf{M}_2=\mathbf{R}_0$.

   \item (3). For a given test level $\alpha=0.05$, the rejection region of test based on the statistic $T$ for testing (\ref{eq6.1})
 is
  $$\{\mathbf{y}_1,\ldots,\mathbf{y}_n: T>t_{\alpha}\},$$
where the critical value $t_{\alpha}$ is obtained by
$$
\alpha=1-\displaystyle\int_{-t_{\alpha}}^{t_{\alpha}}\int_{-t_{\alpha}}^{t_{\alpha}}f(x_1,x_2)\mathrm{d}x_1\mathrm{d}x_2,
$$
with $f(x_1,x_2)$ being the density function of $N\left(\mathbf{0}_2,\left(\begin{array}{ll}1 & \lambda \\ \lambda & 1\end{array}\right)\right).$\\
Here
\begin{equation*}
   \lambda=\dfrac{\sigma_{12}}{\sigma_1\sigma_2}=\dfrac{\sigma_{121}}{\sigma_1\sigma_2}-\dfrac{2\sigma_{122}}{\sigma_1\sigma_2}-\dfrac{2\sigma_{123}}{\sigma_1\sigma_2}+\dfrac{4\sigma_{124}}{\sigma_1\sigma_2},
\end{equation*}
where $\sigma_{121}$,$\sigma_{122}$,$\sigma_{123}$,$\sigma_{124}$ are defined in (\ref{eq5.3}) satisfying that $(\mathbf{M}_1,\mathbf{M}_2,g_{{\ell}_1 1}\\,g_{{\ell}_2 2})$  is equal to $(\mathbf{R}_0^{-1},\mathbf{I}_p,x^2,x^2), (\mathbf{R}_0^{-1},\mathbf{R}_0,x^2,x), (\mathbf{R}_0^{-1}, \mathbf{I}_p,x,x^2), (\mathbf{R}_0^{-1},\mathbf{R}_0\\,x,x)$ , respectively.
\end{itemize}
    
\end{theorem}
Simulation study: Some simulation studies are conducted to evaluate the performance of our proposed statistic $T_1$. The dimension is taken as $p=100, 200, 400,$ and $y=p/n$ is taken as $0.1, 0.5, 0.8, 2.$
The samples $\mathbf{y}_1,\mathbf{y}_2,\ldots,\mathbf{y}_n$ are i.i.d. from a p-dimensional elliptical structure
$$
\mathbf{y}_{j}=\rho_{j}\boldsymbol{\Gamma}\mathbf{x}_j,~j=1,\ldots, n
$$
in (\ref{eq2.1}) where $(\rho_j)^2\sim \chi_p^2$ or $\rho_j\sim Gamma(p,1)$. The test level was set as $\alpha=5\%$, the simulation times were $10000$, and the rejection region  of $T_1$ is $\{\mathbf{y}_1,\ldots,\mathbf{y}_n: \sigma_1^{-1}|T_1-\mu_1|>z_{0.975}\}$, where $z_{0.975}$ is the quantile of $N(0,1)$.

\textbf{Model 1:}
$$
\boldsymbol{\Gamma}=\boldsymbol{\Sigma}^{1/2}, \boldsymbol{\Sigma}=\mathbf{I}_p.
$$

\textbf{Model 2:}
$$
\boldsymbol{\Gamma}=\boldsymbol{\Sigma}^{1/2}, \boldsymbol{\Sigma}=\mathbf{U}(\mathbf{I}_p+\mathbf{D})\mathbf{U}^{\top},
$$
where $\mathbf{U}$ is the eigenvector matrix of $\mathbf{Z}^{\top}\mathbf{Z}$ with the elements $z_{ij}$ being independent and identically distributed (i.i.d.) from~$ N(0,1)
$ and $\mathbf{D}=\textup{diag}(d_{11},d_{22},\ldots,d_{pp})$ being i.i.d. from the uniform distribution $U(0,1)$.

\begin{figure}[htbp]
    \centering
    \begin{subfigure}[b]{0.24\textwidth}
        \centering
        \includegraphics[width=\linewidth]{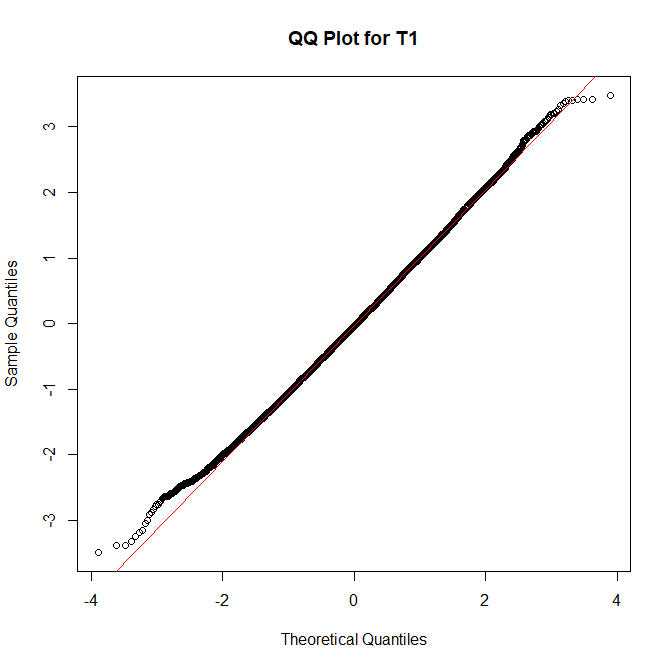}
        \caption{$(p,n)=(100,200)$}
        \label{fig:image1}
    \end{subfigure}%
    \begin{subfigure}[b]{0.24\textwidth}
        \centering
        \includegraphics[width=\linewidth]{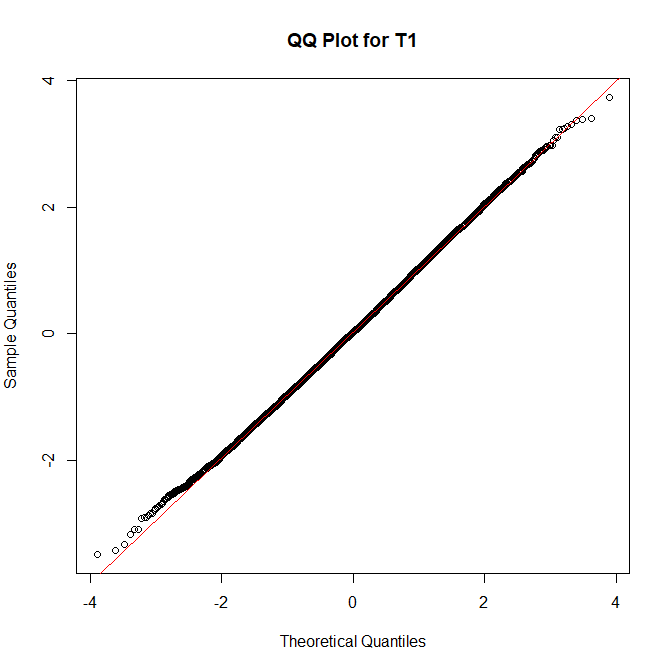}
        \caption{ $(p,n)=(100,200)$}
        \label{fig:image2}
    \end{subfigure}%
    \begin{subfigure}[b]{0.24\textwidth}
        \centering
        \includegraphics[width=\linewidth]{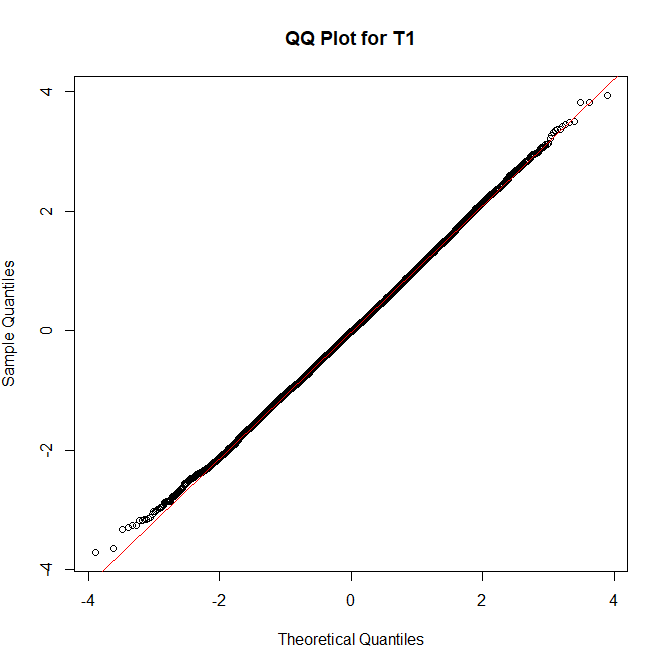}
        \caption{$(p,n)=(250,125)$}
        \label{fig:image3}
    \end{subfigure}%
    \begin{subfigure}[b]{0.24\textwidth}
        \centering
        \includegraphics[width=\linewidth]{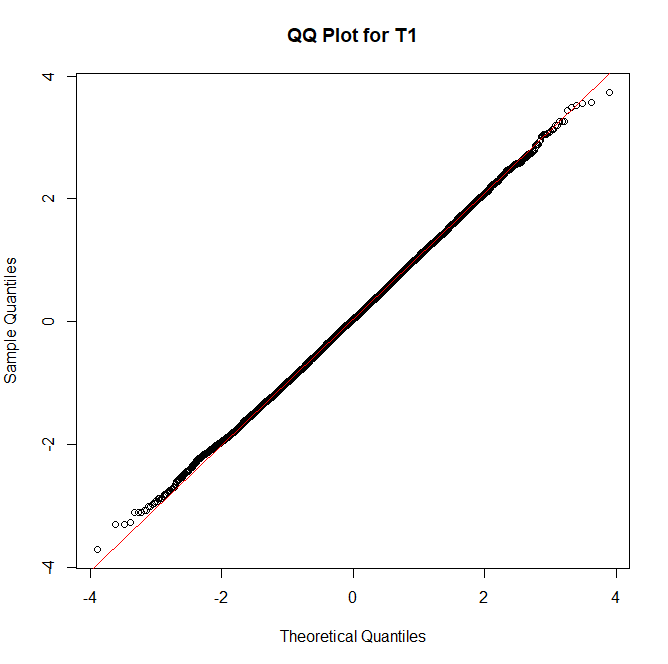}
        \caption{$(p,n)=(250,125)$}
        \label{fig:image4}
    \end{subfigure}
    \caption{Normal QQ-plots for $T_1$ from 10000 independent replications.\\ Note: $\rho\sim Gamma(p,1)$ for (a) and (c); $\rho^2\sim \chi^{2}(p)$ for (b) and (d).}
    \label{fig:four_images}
\end{figure}

Table \ref{table1} present the empirical sizes of $T_1$ for Gaussian population and double exponential population for Model $1$. Table \ref{table2} present the empirical sizes of $T_1$ for Gaussian population and double exponential population for Model $2$. Fig \ref{fig:four_images}  shows the normal QQ-plots of $T_1$ for $(p,n) = (100, 200)$ and $(200, 100) $ for Model $2$.

\begin{table}[htbp]
\normalsize
\begin{center}
\caption{\normalsize Empirical sizes (percentages) of the test $T_1$ for Gaussian and double exponential populations under Model $1$\label{table1}}
\renewcommand{\arraystretch}{1}
\setlength{\tabcolsep}{6pt}
\begin{tabular}{cccc}\hline
$p$~~~~~~~~&$p/n$&~~~~~~~~~Normal  & Double exponential\\ \hline
\multirow{4}{*}{100}~~~~~~~~&$0.1$&~~~~~~~~~4.57&5.68\\
&$0.5$&~~~~~~~~~4.93&5.43\\
&$0.8$&~~~~~~~~~4.86&5.63\\
&$2$&~~~~~~~~~5.3&7.56\\ \hline
\multirow{4}{*}{200}~~~~~~~~&0.1&~~~~~~~~~5.24&5.00\\
&0.5&~~~~~~~~~4.69&4.86\\
&0.8&~~~~~~~~~4.71&5.29\\
&2&~~~~~~~~~5.17&6.63\\ \hline
\multirow{4}{*}{400}~~~~~~~~~&0.1&~~~~~~~~4.83&4.86\\  &0.5&~~~~~~~~~5.01&5.47\\
&0.8&~~~~~~~~~5.16&5.67\\
&2&~~~~~~~~~5.08&5.44\\
\hline
\end{tabular}
\end{center}
\end{table}

\begin{table}[htbp]
\normalsize
\begin{center}
\caption{\normalsize Empirical sizes (percentages) of the test $T_1$ for Gaussian and double exponential populations under Model $2$\label{table2}}
\renewcommand{\arraystretch}{1}
\setlength{\tabcolsep}{6pt}
\begin{tabular}{cccc}\hline
$p$~~~~~~~~&$p/n$&~~~~~~~~~Normal  & Double exponential\\ \hline
\multirow{4}{*}{100}~~~~~~~~&$0.1$&~~~~~~~~~4.78&5.91\\
&$0.5$&~~~~~~~~~4.81&5.61\\
&$0.8$&~~~~~~~~~5.11&5.74\\
&$2$&~~~~~~~~~5.46&6.90\\ \hline
\multirow{4}{*}{200}~~~~~~~~&0.1&~~~~~~~~~5.38&5.17\\
&0.5&~~~~~~~~~4.74&5.04\\
&0.8&~~~~~~~~~4.98&5.40\\
&2&~~~~~~~~~5.29&6.68\\ \hline
\multirow{4}{*}{400}~~~~~~~~~&0.1&~~~~~~~~4.84&5.10\\  &0.5&~~~~~~~~~4.92&5.49\\
&0.8&~~~~~~~~~5.26&5.61\\
&2&~~~~~~~~~4.94&5.29\\
\hline
\end{tabular}
\end{center}
\end{table}
\newpage 

\section{Proofs of Theorem \ref{th2.1}, Theorem \ref{th3.1}, Theorem \ref{th4.1}, Theorem \ref{th3.2} and Example \ref{ex3.1}}\label{appA}

This section provides the skeleton proof of Theorem \ref{th2.1}, Theorem \ref{th3.1} , Theorem \ref{th4.1}, Theorem \ref{th3.2} and Example \ref{ex3.1}. The technical details can be found in the supplementary file: Supplement on ``Spectral properties of high dimensional rescaled sample correlation matrices'.
\subsection{Some preparatory work}
First, notice that
$$
\mathbf{S}_{n}=(n-1)^{-1}\sum_{j=1}^{n}(\mathbf{y}_{j}-\overline{\mathbf{y}})(\mathbf{y}_{j}-\overline{\mathbf{y}})^{\top}=(n-1)^{-1}\sum_{j=1}^{n}(\mathbf{y}_{j}^{0}-\overline{\mathbf{y}}^{0})(\mathbf{y}_{j}^{0}-\overline{\mathbf{y}}^{0})^{\top},
$$
where
\begin{equation*}
\mathbf{y}_{j}^{0}= \left\{
\begin{aligned}
&\rho_{j}\boldsymbol{\Gamma}\mathbf{x}_{j},~~\text{under elliptical case},\\
&\boldsymbol{\Gamma}\mathbf{x}_{j},~~~~~~\text{under linear case},
\end{aligned}
\right.
\end{equation*}
and $\overline{\mathbf{y}}^{0}=n^{-1}\displaystyle\sum_{j=1}^{n}\mathbf{y}_{j}^{0}$. Denote $\mathbf{S}_{n}^{0}=n^{-1}\displaystyle\sum_{j=1}^{n}\mathbf{y}_{j}^{0}{\mathbf{y}_{j}^{0}}^{\top}.$ Since $\mathbf{M}$ is a positive matrix, then based on Theorem A.43 in \cite{bai2010spectral}, we have the LSDs of $\mathbf{S}_n\mathbf{M}$ and $\mathbf{S}_n^0\mathbf{M}$ are the same as $n$ tends to infinity. Therefore, we consider $\mathbf{S}_{n}^{0}$ instead of $\mathbf{S}_{n}$ in the proof of Theorem \ref{th2.1} in the next section. Similar to the discussion in Section A.1 of \cite{yin2023central}, combined with the positive definiteness of $\mathbf{M}$, the substitution principle for the CLT of the LSS of $\widehat{\mathbf{R}}_{n}\mathbf{M}$ holds. That is, we could study the matrix $\mathbf{S}_{n}^{0}$ instead of $\mathbf{S}_{n}$ in the proof of Theorems \ref{th3.1} and \ref{th4.1}.


In the following sections, we still use $\mathbf{S}_{n}$ instead of $\mathbf{S}_{n}^{0}$, and the definition of $\widehat{\mathbf{R}}_{n}\mathbf{M}$ is redefined accordingly. Before presenting the proof, we give some notations. Let
\begin{equation}
\mathbf{G}=(g_{ki})=\left[\textup{diag}(\boldsymbol{\Sigma})\right]^{-1/2}\boldsymbol{\Gamma},
\end{equation}
\begin{equation}
\boldsymbol{\Xi}_{j}=\begin{cases}
\rho_{j}^{2}\mathbf{G}\mathbf{x}_{j}\mathbf{x}_{j}^{\top}\mathbf{G}^{\top}, &\text{under elliptical case},\\
\mathbf{G}\mathbf{x}_{j}\mathbf{x}_{j}^{\top}\mathbf{G}^{\top}, &\text{under linear case},
\end{cases}
\end{equation}
\begin{equation}
\boldsymbol{\Xi}=\dfrac{1}{n}\displaystyle\sum_{j=1}^{n}\boldsymbol{\Xi}_{j}.
\end{equation}
Observe that
\begin{equation}
\textup{E}(\boldsymbol{\Xi}_{j})=\textup{E}(\boldsymbol{\Xi})=\mathbf{R},~~~\text{hence} ~\textup{diag}(\textup{E}(\boldsymbol{\Xi}_{j}))=\textup{diag}(\textup{E}(\boldsymbol{\Xi}))=\mathbf{I}_{p}.
\end{equation}
It is easy to see that $\widehat{\mathbf{R}}_{n}\mathbf{M}$ can also be written as
\begin{equation*}
\widehat{\mathbf{R}}_{n}\mathbf{M}=\left[\textup{diag}(\boldsymbol{\Xi})\right]^{-1/2}\boldsymbol{\Xi}\left[\textup{diag}(\boldsymbol{\Xi})\right]^{-1/2}\mathbf{M}.
\end{equation*}
Throughout this paper, we also note that $C$ and $C_{(\cdot)}$ denote constants that may take different values from one appearance to another.

\subsection{Proof of Theorem \ref{th2.1}}

First, in the elliptical case, from \cite{yin2023central} and Assumption B, we can easily obtain that
\begin{equation}
\Vert \widehat{\mathbf{R}}_{n}\mathbf{M}-\boldsymbol{\Xi}\mathbf{M}\Vert\leq \Vert \widehat{\mathbf{R}}_{n}-\boldsymbol{\Xi}\Vert\cdot\Vert\mathbf{M}\Vert\to 0. ~~~a.s..
\end{equation}
Therefore, we need to focus only on the spectrum of the matrix $\boldsymbol{\Xi}\mathbf{M}$ or $\mathbf{M}^{1/2}\boldsymbol{\Xi}\mathbf{M}^{1/2}$.

\noindent Then, the result follows from Theorem 2.1 in \cite{yin2023central} and Weyl's inequality.

Now, we consider the linear case, according to Lemma 4 of \cite{el2009concentration} and Assumption B, we have
\begin{align*}
   &\Vert \widehat{\mathbf{R}}_{n}\mathbf{M}-\boldsymbol{\Xi}\mathbf{M}\Vert\\
   \leq &\Vert \widehat{\mathbf{R}}_{n}-\boldsymbol{\Xi}\Vert\cdot\Vert\mathbf{M}\Vert\\
   \leq &\Vert(\textup{diag}(\boldsymbol{\Xi}))^{-1/2}-\mathbf{I}_{p}\Vert^{2}\cdot\Vert\boldsymbol{\Xi}\Vert\cdot\Vert\mathbf{M}\Vert+2\Vert(\textup{diag}(\boldsymbol{\Xi}))^{-1/2}-\mathbf{I}_{p}\Vert\cdot\Vert\boldsymbol{\Xi}\Vert\cdot\Vert\mathbf{M}\Vert\\
   \to& 0~~~a.s..
\end{align*}
Then, the result of this theorem in the linear case follows from Theorem 1 in \cite{el2009concentration} and Weyl's inequality.

\subsection{Sketch of proofs of Theorem \ref{th3.1} and Theorem \ref{th4.1}}

This section provides the main sketch of the proof of Theorems \ref{th3.1} and \ref{th4.1}. The details are included in the Supplementary Material. Recall that, for any analaytic $g$ in a domain containing the support interval $F^{y,H}$,
\begin{equation*}
    W(g)=\displaystyle\sum_{i=1}^{p}g(\widehat{\lambda}_{i})-p\displaystyle\int g(x)\textup{d}F^{y_{n},H_{n}},
\end{equation*}
where $\{ \widehat{\lambda}_{i}, i=1,\ldots, p\}$ are the eigenvalues of $\widehat{\mathbf{R}}_{n}\mathbf{M}$.

\subsubsection{Truncation, centralization and rescaling}

First, in the elliptical case, we begin the proof of Theorem \ref{th3.1} by truncating the variable $\rho$ at a proper order of $n$. From the moment condition $\textup{E}\left\vert\dfrac{\rho_{j}^{2}-p}{\sqrt{p}} \right\vert ^{2+\varepsilon}<\infty$ for some $\varepsilon > 0$ in Assumption $A_{E}$, we can choose a sequence of $\eta_{n}\downarrow  0$ such that
\begin{equation*}
       \eta_{n}\sqrt{n}\to \infty,~\eta_{n}^{-2}p^{-1}\textup{E}\left[(\rho_{1}^{2}-p)^{2}I_{\{|\rho_{1}^{2}-p|\geq \eta_{n}p\}}\right] \to 0.
\end{equation*}
Denote $\check{\boldsymbol{\Xi}}= n^{-1} \displaystyle\sum_{j=1}^{n}\check{\rho}_{j}^{2}\mathbf{G}\mathbf{x}_{j}\mathbf{x}_{j}^{\top}\mathbf{G}^{\top}$ where $\check{\rho}_{j}=\rho_{j} I_{\{|\rho_{j}^{2}-p|<\eta_{n}p\}}$,
\begin{equation*}
    \check{\mathbf{R}}_{n}\mathbf{M}=\left[ \textup{diag}(\check{\boldsymbol{\Xi}})\right]^{-1/2}\check{\boldsymbol{\Xi}}\left[ \textup{diag}(\check{\boldsymbol{\Xi}})\right]^{-1/2}\mathbf{M}
\end{equation*}
and $\check{W}(g)$ is the truncated version of $W(g)$. By \cite{yin2023central} we have
\begin{equation*}
    \textup{P}(\widehat{\mathbf{R}}_{n}\mathbf{M}\neq \check{\mathbf{R}}_{n}\mathbf{M}, i.o.)=\textup{P}(\widehat{\mathbf{R}}_{n}\neq \check{\mathbf{R}}_{n}, i.o.)\to 0.
\end{equation*}
Now define $\widetilde{\boldsymbol{\Xi}}=n^{-1} \displaystyle\sum_{j=1}^{n}\widetilde{\rho}_{j}^{2}\mathbf{G}\mathbf{x}_{j}\mathbf{x}_{j}^{\top}\mathbf{G}^{\top}$ where $\widetilde{\rho_{j}}=\dfrac{\check{\rho}_{j}}{\sigma_{n}}$ and $\sigma_{n}^{2}=\dfrac{\textup{E}(\check{\rho}_{1}^{2})}{p}$. Also define $\widetilde{\mathbf{R}}_{n}\mathbf{M}$ and $\widetilde{W}(f)$ as the analogues of  $\mathbf{R}_{n}\mathbf{M}$ and $W(f)$ with $\rho_{j}$ being replaced by $\widetilde{\rho}_{j}$. By \cite{yin2023central} we get
\begin{equation*}
    \left\Vert \check{\mathbf{R}}_{n}\mathbf{M}-\widetilde{\mathbf{R}}_{n}\mathbf{M}\right\Vert\leq \left\Vert \check{\mathbf{R}}_{n}-\widetilde{\mathbf{R}}_{n}\right\Vert \cdot \left\Vert \mathbf{M} \right\Vert =o_{a.s.}(n^{-1}).
\end{equation*}
We finally obtain for large $n$
\begin{equation}
   |\check{W}(g)-\widetilde{W}(g)|\leq C \displaystyle\sum_{i=1}^{p}\left\vert\lambda_{i}(\check{\mathbf{R}}_{n}\mathbf{M})-\lambda_{i}(\widetilde{\mathbf{R}}_{n}\mathbf{M})\right\vert\leq Cp\cdot \left\Vert\check{\mathbf{R}}_{n}\mathbf{M}-\widetilde{\mathbf{R}}_{n}\mathbf{M}\right\Vert=o_{a.s.}(1) ,
\end{equation}
where $C$ is a bound on $|g'(z)|$ . Therefore, we only need to find the limiting distribution of $\widetilde{W}(g)$. For simplicity, we still use $W(g), \widehat{\mathbf{R}}_{n}\mathbf{M}, \rho_{j}$ instead of $\widetilde{W}(g), \widetilde{\mathbf{R}}_{n}\mathbf{M}, \widetilde{\rho}_{j}$ respectively. And assume that
\begin{equation}
    \forall j, |\rho_{j}^{2}-p|<\eta_{n}p,~~\textup{E}(\rho_{j}^{2})=p,~~\textup{E}(\rho_{j}^{4})=p^{2}+\tau p +o(p),
\end{equation}
in the subsequent discussion.

Then, in the linear case, we perform truncation, centralization, and rescaling on $\{x_{ij},i=1,\ldots,p, j=1,\ldots, n\}$ and provide the following notations.

Denote ~$\check{\mathbf{X}}=(\check{x}_{ij})$~with $\check{x}_{ij}=x_{ij}I_{\{|x_{ij}|<\eta_{n}\sqrt{n}\}},~\check{\boldsymbol{\Xi}}=n^{-1}\mathbf{G}\check{\mathbf{X}}\check{\mathbf{X}}^{\top}\mathbf{G}^{\top},$
\begin{equation*}
    \check{\mathbf{R}}_{n}\mathbf{M}=n^{-1}\left[\textup{diag}(\check{\boldsymbol{\Xi}})\right]^{-1/2}\mathbf{G}\mathbf{\check{X}}\mathbf{\check{X}}^{\top}\mathbf{G}^{\top}\left[\textup{diag}(\check{\boldsymbol{\Xi}})\right]^{-1/2}\mathbf{M}
\end{equation*}
and $\check{W}(g)$ is the truncated version of $W(g)$. As have been proved in \cite{el2009concentration}, under the moment assumption, we shall select a sequence of $\eta_{n}=(\textup{log}~n)^{-(1+\varepsilon)/2}\to 0$ as $n\to \infty$ satisfying
\begin{equation*}
    \textup{P}(\widehat{\mathbf{R}}_{n}\mathbf{M}\neq \check{\mathbf{R}}_{n}\mathbf{M}, i.o.)=\textup{P}(\widehat{\mathbf{R}}_{n}\neq \check{\mathbf{R}}_{n}, i.o.)\to 0.
\end{equation*}
Now define $\widetilde{\mathbf{X}}=(\tilde{x}_{ij})$ with $\tilde{x}_{ij}=\left(\check{x}_{ij}-\textup{E}~\check{x}_{ij}\right)/\sqrt{\textup{E}\left(\check{x}_{ij}-\textup{E}~\check{x}_{ij}\right)^{2}}$. Also define $\widetilde{\boldsymbol{{\Xi}}},~\widetilde{\mathbf{R}}_{n}\mathbf{M}$ and $\widetilde{W}(f)$ as the analogues of $\boldsymbol{\Xi},~\mathbf{R}_{n}\mathbf{M}$ and $W(f)$ with $\mathbf{X}$ being replaced by $\widetilde{\mathbf{X}}$. For large $n$ and any $1\leq i\leq p, 1\leq j \leq n,$ by \cite{yin2023central} we get
\begin{equation}
    \left\Vert \check{\mathbf{R}}_{n}\mathbf{M}-\widetilde{\mathbf{R}}_{n}\mathbf{M}\right\Vert\leq \left\Vert \check{\mathbf{R}}_{n}-\widetilde{\mathbf{R}}_{n}\right\Vert \cdot \left\Vert \mathbf{M} \right\Vert =o_{a.s.}(n^{-1}).
\end{equation}
We finally obtain for large n
\begin{equation}
   \left|\check{W}(g)-\widetilde{W}(g)\right|\leq C \displaystyle\sum_{i=1}^{p}\left\vert\lambda_{i}(\check{\mathbf{R}}_{n}\mathbf{M})-\lambda_{i}(\widetilde{\mathbf{R}}_{n}\mathbf{M})\right\vert\leq Cp\cdot \left\Vert\check{\mathbf{R}}_{n}\mathbf{M}-\widetilde{\mathbf{R}}_{n}\mathbf{M}\right\Vert=o_{a.s.}(1) .
\end{equation}
where $C$ is a bound on $|g'(z)|$ .

Therefore, we shall assume in the following that the underlying variables in the data matrix $\mathbf{X}$ are all truncated at $\eta_{n}\sqrt{n}$, centralized and rescalized to have unit variances.

The proofs of Theorems \ref{th3.1} and \ref{th4.1} rely on analyzing the Stieltjes transform  $s_{n}(z)$ of the ESD of $\widehat{\mathbf{R}}_{n}\mathbf{M}$. We denote $M_{n}(z):= p(s_{n}(z)-s_{y_{n}}(z)),$ where $s_{y_{n}}(z)$ is the Stieltjes transform of the distribution $F^{y_{n},H_{n}}$. Notice that by the Cauchy integral formula, we have
\begin{equation}\label{eqA}
W(g)=\displaystyle\sum_{i=1}^{p}g(\widehat{\lambda}_{i})-p\displaystyle\int g(x) \textup{d}F^{y_{n},H_{n}}=-\dfrac{1}{2\pi i}\oint_{\mathcal{C}}g(z)M_{n}(z)\textup{d}z,
\end{equation}
where $\mathcal{C}$ is any contour inside the domain and surrounding the support interval of $F^{y,H}$. This suggests that our target is to analyze the random process $M_{n}(z).$ Following the ideas of the arguments on pages 1000-1001 in \cite{yin2023central}, we investigate a truncated version $\widehat{M}_{n}(z)$ of $M_{n}(z)$. Let $\mathbf{x}_{r}$ be any number greater than $\displaystyle\limsup_{n}\lambda_{\textup{max}}^{\mathbf{RM}}(1+\sqrt{y})^{2},$ and $\mathbf{x}_{\ell}$ be any negative number if $\displaystyle\liminf_{n} \lambda_{\min}^{\mathbf{RM}}I_{(0, 1)}(y)(1-\sqrt{y})^{2}=0.$ Otherwise choose $\mathbf{x}_{\ell}\in (0, \displaystyle\liminf_{n}\lambda_{\min}^{\mathbf{RM}}I_{(0, 1)}(y)(1-\sqrt{y})^{2}).$ Let $v_{0}>0$ be arbitrary, we define a counter $\mathcal{C}$ as $\mathcal{C}\equiv\mathcal{C}_{\ell}\cup\mathcal{C}_{r}\cup\mathcal{C}_{u}\cup\mathcal{C}_{b}$, where
\begin{align*}
    \mathcal{C}_{\ell}&=\{x_{\ell}+iv: |v|\leq v_{0}\},~~~
    \mathcal{C}_{u}=\{x+iv_{0}: x\in[x_{\ell},x_{r}]\}\\
    \mathcal{C}_{r}&=\{x_{r}+iv: |v|\leq v_{0}\},~~~
    \mathcal{C}_{b}=\{x-iv_{0}: x\in[x_{\ell},x_{r}]\}.
\end{align*}
Then, we define the subsets $\mathcal{C}_{n}$ of $\mathcal{C}$ as follows
\begin{equation*}
    \mathcal{C}_{n}=\mathcal{C}\cap\{z: |\mathfrak{I}z|>n^{-2}\}.
\end{equation*}
For $z=x+iv$, we define
\begin{equation*}
\widehat{M}_{n}(z)= \left\{
\begin{aligned}
&M_{n}(z),~~~~~~~~~~~~~~~~~~~~\text{if}~z\in\mathcal{C}_{n},\\
&M_{n}(x_{\ell}+in^{-2}),~~~~~~~\text{if}~x=x_{\ell}, v\in[0, n^{-2}],\\
&M_{n}(x_{\ell}-in^{-2}),~~~~~~~\text{if}~x=x_{\ell}, v\in[-n^{-2}, 0),\\
&M_{n}(x_{r}+in^{-2}),~~~~~~~\text{if}~x=x_{r}, v\in[0, n^{-2}],\\
&M_{n}(x_{r}-in^{-2}),~~~~~~~\text{if}~x=x_{r}, v\in[-n^{-2}, 0).
\end{aligned}
\right.
\end{equation*}
Note that $\widehat{M}_{n}(z)$ agrees with $M_{n}(z)$ on $\mathcal{C}_{n}$. Similar to the discussion in \cite{yin2023central}, we have
\begin{equation*}
    \left\vert \displaystyle\oint_{\mathcal{C}}g(z)(M_{n}(z)-\widehat{M}_{n}(z))\textup{d}z \right\vert=\left\vert \displaystyle\oint_{\mathcal{C}-\mathcal{C}_{n}}g(z)(M_{n}(z)-\widehat{M}_{n}(z))\textup{d}z \right\vert=o(p^{-1}).
\end{equation*}
Therefore, we only need to focus on $\displaystyle\oint_{\mathcal{C}}g(z)\widehat{M}_{n}(z)\textup{d}z.$ Since $v_{0}$ can be chosen arbitrarily small, the contributions from segments $\mathcal{C}_{\ell}$ and $\mathcal{C}_{r}$ can also be small. As a result, we only need to focus on $z\in \mathcal{C}_{u}\cup \mathcal{C}_{b}$ when analyzing $\widehat{M}_{n}(z)$. For simplicity, we still use $M_{n}(z)$ instead of $\widehat{M}_{n}(z)$ in the following sections.

\subsubsection{Analysis of the random process $M_n(z)$}

We split $M_{n}(z)$ into several parts. Specifically, we have
\begin{equation}
M_{n}(z)=V+M_{0}+M_{1}+M_{2}+O,
\end{equation}
where
\begin{align*}
V&=\Big[\textup{tr}\mathbf{A}^{-1}(z)-\textup{Etr}\mathbf{A}^{-1}(z)\Big]+\dfrac{1}{2}\Big[\textup{tr}(\mathbf{A}^{-1}(z)\mathbf{D})-\textup{Etr}(\mathbf{A}^{-1}(z)\mathbf{D})\Big]\\
&+\dfrac{z}{2}\Big[\textup{tr}(\mathbf{A}^{-2}(z)\mathbf{D})-\textup{Etr}(\mathbf{A}^{-2}(z)\mathbf{D})\Big]\\
&+\dfrac{1}{2}\Big[\textup{tr}(\mathbf{MA}^{-1}(z)\mathbf{M}^{-1}\mathbf{D})-\textup{Etr}(\mathbf{MA}^{-1}(z)\mathbf{M}^{-1}\mathbf{D})\Big]\\
 &+\dfrac{z}{2}\Big[\textup{tr}(\mathbf{M}\mathbf{A}^{-2}(z)\mathbf{M}^{-1}\mathbf{D})-\textup{Etr}(\mathbf{MA}^{-2}(z)\mathbf{M}^{-1}\mathbf{D})\Big],\\
M_{1}&=-\dfrac{1}{8}\textup{tr}(\mathbf{A}^{-1}(z)\mathbf{D}^{2})-\dfrac{1}{8}\textup{tr}(\mathbf{A}^{-1}(z)\mathbf{M}^{-1}\mathbf{D}^{2}\mathbf{M})-\dfrac{z}{8}\textup{tr}(\mathbf{A}^{-2}(z)\mathbf{D}^{2})\\
&-\dfrac{z}{8}\textup{tr}(\mathbf{A}^{-2}(z)\mathbf{M}^{-1}\mathbf{D}^{2}\mathbf{M})+\dfrac{z}{4}\textup{tr}(\mathbf{A}^{-1}(z)\mathbf{D})^{2}+\dfrac{z}{2}\textup{tr}(\mathbf{A}^{-1}(z)\mathbf{M}^{-1}\mathbf{D}\mathbf{M}\mathbf{A}^{-1}(z)\mathbf{D})\\
&+\dfrac{z}{4}\textup{tr}(\mathbf{A}^{-1}(z)\mathbf{M}^{-1}\mathbf{D}\mathbf{M})^{2}+\dfrac{1}{4}\textup{tr}(\mathbf{A}^{-1}(z)\mathbf{D}\mathbf{M}^{-1}\mathbf{D}\mathbf{M})+\dfrac{z}{4}\textup{tr}(\mathbf{A}^{-2}(z)\mathbf{D}\mathbf{M}^{-1}\mathbf{D}\mathbf{M}),\\
M_{2}&=\dfrac{z^{2}}{4}\textup{tr}\left[\mathbf{A}^{-1}(z)(\mathbf{A}^{-1}(z)\mathbf{D})^{2}\right]+\dfrac{z^{2}}{4}\textup{tr}\left[\mathbf{A}^{-1}(z)\mathbf{D}\mathbf{A}^{-2}(z)\mathbf{M}^{-1}\mathbf{D}\mathbf{M}\right]\\
&+\dfrac{z^{2}}{4}\textup{tr}\left[\mathbf{A}^{-2}(z)\mathbf{D}\mathbf{A}^{-1}(z)\mathbf{M}^{-1}\mathbf{D}\mathbf{M}\right]+\dfrac{z^{2}}{4}\textup{tr}\left[\mathbf{A}^{-1}(z)(\mathbf{A}^{-1}(z)\mathbf{M}^{-1}\mathbf{D}\mathbf{M})^{2}\right],\\
M_{0}&=\textup{E}\textup{tr}\mathbf{A}^{-1}(z)-ps_{y_{n}}(z)+\dfrac{1}{2}\textup{E}\textup{tr}(\mathbf{A}^{-1}(z)\mathbf{D})+\dfrac{z}{2}\textup{E}\textup{tr}(\mathbf{A}^{-2}(z)\mathbf{D})\\
&+\dfrac{1}{2}\textup{E}\textup{tr}(\mathbf{M}\mathbf{A}^{-1}(z)\mathbf{M}^{-1}\mathbf{D})+\dfrac{z}{2}\textup{E}\textup{tr}(\mathbf{M}\mathbf{A}^{-2}(z)\mathbf{M}^{-1}\mathbf{D}),
\end{align*}
O is the trace of all terms that contain $\mathbf{D}$ several times greater than or equal to 3, and $\mathbf{A}(z)=\boldsymbol{\Xi}\mathbf{M}-z\mathbf{I}_{p}, \mathbf{D}=\textup{diag}(\boldsymbol{\Xi})-\mathbf{I}_{p}$. Consequently,

(1): In both elliptical and linear cases, the term O converges in probability to zero; thus, it has no contribution to the limit properties of $M_{n}(z)$.

(2): In both elliptical and linear cases, the terms $M_{1}$ and $M_{2}$ converge in probability to the limit of their expectation, but do not contribute to the limit of the variance-covariance function of $M_{n}(z)$. Moreover, we have in the elliptical case
\begin{align*}
&\textup{E}M_{1}+\textup{E}M_{2}\to\textup{E}_{1}(z)\\
=&\dfrac{1}{4z}\displaystyle\lim_{n\to \infty}\dfrac{1}{n}\sum_{\ell=1}^{p}\left[\mathbf{e}_{\ell}^{\top}\left(\mathbf{I}+\underline{s}(z)\mathbf{R}\mathbf{M}\right)^{-1}\mathbf{e}_{\ell}\right]+\dfrac{1}{4z}\displaystyle\lim_{n\to \infty}\dfrac{1}{n}\sum_{\ell=1}^{p}\left[\mathbf{e}_{\ell}^{\top}\mathbf{M}\left(\mathbf{I}+\underline{s}(z)\mathbf{R}\mathbf{M}\right)^{-1}\mathbf{M}^{-1}\mathbf{e}_{\ell}\right]\\
-&\dfrac{1}{4z}\displaystyle\lim_{n\to \infty}\dfrac{1}{n}\sum_{\ell=1}^{p}\left[\mathbf{e}_{\ell}^{\top}\left(\mathbf{I}+\underline{s}(z)\mathbf{R}\mathbf{M}\right)^{-2}\mathbf{e}_{\ell}\right]-\dfrac{1}{4z}\displaystyle\lim_{n\to \infty}\dfrac{1}{n}\sum_{\ell=1}^{p}\left[\mathbf{e}_{\ell}^{\top}\mathbf{M}\left(\mathbf{I}+\underline{s}(z)\mathbf{R}\mathbf{M}\right)^{-2}\mathbf{M}^{-1}\mathbf{e}_{\ell}\right]\\
-&\dfrac{1}{2z}\displaystyle\lim_{n\to \infty}\dfrac{1}{n}\sum_{k,\ell=1}^{p}r_{k\ell}^{2}\mathbf{e}_{\ell}^{\top}\mathbf{M}^{-1}\mathbf{e}_{k}\left[\mathbf{e}_{k}^{\top}\mathbf{M}\left(\mathbf{I}+\underline{s}(z)\mathbf{R}\mathbf{M}\right)^{-1}\mathbf{e}_{\ell}\right]\\
+&\dfrac{1}{2z}\displaystyle\lim_{n\to \infty}\dfrac{1}{n}\sum_{k,\ell=1}^{p}r_{k\ell}^{2}\mathbf{e}_{\ell}^{\top}\mathbf{M}^{-1}\mathbf{e}_{k}\left[\mathbf{e}_{k}^{\top}\mathbf{M}\left(\mathbf{I}+\underline{s}(z)\mathbf{R}\mathbf{M}\right)^{-2}\mathbf{e}_{\ell}\right]\\
+&\dfrac{1}{4}\displaystyle\lim_{n\to \infty}\dfrac{1}{n}\sum_{k,\ell=1}^{p}r_{k\ell}^{2}\dfrac{\partial}{\partial z}\left[\mathbf{e}_{k}^{\top}\left(\mathbf{I}+\underline{s}(z)\mathbf{R}\mathbf{M}\right)^{-1}\mathbf{e}_{\ell}\cdot \mathbf{e}_{\ell}^{\top}\left(\mathbf{I}+\underline{s}(z)\mathbf{R}\mathbf{M}\right)^{-1}\mathbf{e}_{k}\right]\\
+&\dfrac{1}{2}\displaystyle\lim_{n\to \infty}\dfrac{1}{n}\sum_{k,\ell=1}^{p}r_{k\ell}^{2}\dfrac{\partial}{\partial z}\left[\mathbf{e}_{k}^{\top}\left(\mathbf{I}+\underline{s}(z)\mathbf{R}\mathbf{M}\right)^{-1}\mathbf{M}^{-1}\mathbf{e}_{\ell}\cdot \mathbf{e}_{\ell}^{\top}\mathbf{M}\left(\mathbf{I}+\underline{s}(z)\mathbf{R}\mathbf{M}\right)^{-1}\mathbf{e}_{k}\right]\\
+&\dfrac{1}{4}\displaystyle\lim_{n\to \infty}\dfrac{1}{n}\sum_{k,\ell=1}^{p}r_{k\ell}^{2}\dfrac{\partial}{\partial z}\left[\mathbf{e}_{k}^{\top}\mathbf{M}\left(\mathbf{I}+\underline{s}(z)\mathbf{R}\mathbf{M}\right)^{-1}\mathbf{M}^{-1}\mathbf{e}_{\ell}\cdot \mathbf{e}_{\ell}^{\top}\mathbf{M}\left(\mathbf{I}+\underline{s}(z)\mathbf{R}\mathbf{M}\right)^{-1}\mathbf{M}^{-1}\mathbf{e}_{k}\right].
\end{align*}
While in the linear case,
\begin{align*}
&\textup{E}M_{1}+\textup{E}M_{2}\to\textup{E}_{1}(z)\\
=&\dfrac{1}{4z}\displaystyle\lim_{n\to \infty}\dfrac{1}{n}\sum_{\ell=1}^{p}\left[\mathbf{e}_{\ell}^{\top}\left(\mathbf{I}+\underline{s}(z)\mathbf{R}\mathbf{M}\right)^{-1}\mathbf{e}_{\ell}\right]\cdot\left[\dfrac{\beta_{x}}{2}\displaystyle\sum_{j=1}^{p}\left(\mathbf{e}_{\ell}^{\top}\mathbf{G}\mathbf{e}_{j}\right)^{4}+1\right]\\
+&\dfrac{1}{4z}\displaystyle\lim_{n\to \infty}\dfrac{1}{n}\sum_{\ell=1}^{p}\left[\mathbf{e}_{\ell}^{\top}\mathbf{M}\left(\mathbf{I}+\underline{s}(z)\mathbf{R}\mathbf{M}\right)^{-1}\mathbf{M}^{-1}\mathbf{e}_{\ell}\right]\cdot\left[\dfrac{\beta_{x}}{2}\displaystyle\sum_{j=1}^{p}\left(\mathbf{e}_{\ell}^{\top}\mathbf{G}\mathbf{e}_{j}\right)^{4}+1\right]\\
-&\dfrac{1}{4z}\displaystyle\lim_{n\to \infty}\dfrac{1}{n}\sum_{\ell=1}^{p}\left[\mathbf{e}_{\ell}^{\top}\left(\mathbf{I}+\underline{s}(z)\mathbf{R}\mathbf{M}\right)^{-2}\mathbf{e}_{\ell}\right]\cdot\left[\dfrac{\beta_{x}}{2}\displaystyle\sum_{j=1}^{p}\left(\mathbf{e}_{\ell}^{\top}\mathbf{G}\mathbf{e}_{j}\right)^{4}+1\right]\\
-&\dfrac{1}{4z}\displaystyle\lim_{n\to \infty}\dfrac{1}{n}\sum_{\ell=1}^{p}\left[\mathbf{e}_{\ell}^{\top}\mathbf{M}\left(\mathbf{I}+\underline{s}(z)\mathbf{R}\mathbf{M}\right)^{-2}\mathbf{M}^{-1}\mathbf{e}_{\ell}\right]\cdot\left[\dfrac{\beta_{x}}{2}\displaystyle\sum_{j=1}^{p}\left(\mathbf{e}_{\ell}^{\top}\mathbf{G}\mathbf{e}_{j}\right)^{4}+1\right]\\
-&\dfrac{1}{2z}\displaystyle\lim_{n\to \infty}\dfrac{1}{n}\displaystyle\sum_{k,\ell=1}^{p}\mathbf{e}_{\ell}^{\top}\mathbf{M}^{-1}\mathbf{e}_{k}\left[\mathbf{e}_{k}^{\top}\mathbf{M}\left(\mathbf{I}+\underline{s}(z)\mathbf{R}\mathbf{M}\right)^{-1}\mathbf{e}_{\ell}\right]\cdot \left[r_{k\ell}^{2}+\dfrac{\beta_{x}}{2}\displaystyle\sum_{j=1}^{p}\left(\mathbf{e}_{\ell}^{\top}\mathbf{G}\mathbf{e}_{j}\mathbf{e}_{j}^{\top}\mathbf{G}^{\top}\mathbf{e}_{k}\right)^{2}\right]\\
+&\dfrac{1}{2z}\displaystyle\lim_{n\to \infty}\dfrac{1}{n}\displaystyle\sum_{k,\ell=1}^{p}\mathbf{e}_{\ell}^{\top}\mathbf{M}^{-1}\mathbf{e}_{k}\left[\mathbf{e}_{k}^{\top}\mathbf{M}\left(\mathbf{I}+\underline{s}(z)\mathbf{R}\mathbf{M}\right)^{-2}\mathbf{e}_{\ell}\right]\cdot \left[r_{k\ell}^{2}+\dfrac{\beta_{x}}{2}\displaystyle\sum_{j=1}^{p}\left(\mathbf{e}_{\ell}^{\top}\mathbf{G}\mathbf{e}_{j}\mathbf{e}_{j}^{\top}\mathbf{G}^{\top}\mathbf{e}_{k}\right)^{2}\right]\\
+&\dfrac{1}{4}\lim_{n\to \infty}\dfrac{1}{n}\displaystyle\sum_{k,\ell=1}^{p}\dfrac{\partial}{\partial z}\bigg[2\mathbf{e}_{k}^{\top}\left(\mathbf{I}+\underline{s}(z)\mathbf{R}\mathbf{M}\right)^{-1}\mathbf{M}^{-1}\mathbf{e}_{\ell}\mathbf{e}_{\ell}^{\top}\mathbf{M}\left(\mathbf{I}+\underline{s}(z)\mathbf{R}\mathbf{M}\right)^{-1}\mathbf{e}_{k}\\
+&\mathbf{e}_{k}^{\top}\mathbf{M}\left(\mathbf{I}+\underline{s}(z)\mathbf{R}\mathbf{M}\right)^{-1}\mathbf{M}^{-1}\mathbf{e}_{\ell}\mathbf{e}_{\ell}^{\top}\mathbf{M}\left(\mathbf{I}+\underline{s}(z)\mathbf{R}\mathbf{M}\right)^{-1}\mathbf{M}^{-1}\mathbf{e}_{k}+\mathbf{e}_{k}^{\top}\left(\mathbf{I}+\underline{s}(z)\mathbf{R}\mathbf{M}\right)^{-1}\mathbf{e}_{\ell}\\
&\cdot\mathbf{e}_{\ell}^{\top}\left(\mathbf{I}+\underline{s}(z)\mathbf{R}\mathbf{M}\right)^{-1}\mathbf{e}_{k}\bigg]\cdot \left[r_{k\ell}^{2}+\dfrac{\beta_{x}}{2}\displaystyle\sum_{j=1}^{p}\left(\mathbf{e}_{\ell}^{\top}\mathbf{G}\mathbf{e}_{j}\mathbf{e}_{j}^{\top}\mathbf{G}^{\top}\mathbf{e}_{k}\right)^{2}\right].
\end{align*}

(3): The term $M_{0}$ converges to certain limits in both cases. The limit of $M_{0}$ in the elliptical case is
\begin{align*}
\textup{E}_{0}(z)=&y\displaystyle\int\dfrac{(\underline{s}'(z)t)^{2}\mathrm{d}H(t)}{\underline{s}(z)(1+t\underline{s}(z))^{3}}+(\tau -2)(1+z\underline{s}(z))\displaystyle\int\dfrac{\underline{s}'(z)t\mathrm{d}H(t)}{(1+t\underline{s}(z))^{2}}\\
+&\displaystyle\lim_{n\to \infty } \dfrac{1}{n} \sum_{k=1}^{p}\dfrac{\partial}{\partial z}\left[\underline{s}(z)\mathbf{e}_{k}^{\top}\left(\mathbf{I}+\underline{s}(z)\mathbf{R}\mathbf{M}\right)^{-1}\mathbf{R}\mathbf{e}_{k}\cdot \mathbf{e}_{k}^{\top}\mathbf{R}\mathbf{M}\left(\mathbf{I}+\underline{s}(z)\mathbf{R}\mathbf{M}\right)^{-1}\mathbf{e}_{k}\right]\\
+&\displaystyle\lim_{n\to \infty } \dfrac{1}{n} \sum_{k=1}^{p}\dfrac{\partial}{\partial z}\left[\underline{s}(z)\mathbf{e}_{k}^{\top}\mathbf{M}\left(\mathbf{I}+\underline{s}(z)\mathbf{R}\mathbf{M}\right)^{-1}\mathbf{R}\mathbf{e}_{k}\cdot \mathbf{e}_{k}^{\top}\mathbf{R}\mathbf{M}\left(\mathbf{I}+\underline{s}(z)\mathbf{R}\mathbf{M}\right)^{-1}\mathbf{M}\mathbf{e}_{k}\right].
\end{align*}
The limit of $M_{0}$ in the linear case is
\begin{align*}
\textup{E}_{0}(z)=&y\displaystyle\int\dfrac{(\underline{s}'(z)t)^{2}\mathrm{d}H(t)}{\underline{s}(z)(1+t\underline{s}(z))^{3}}\\
+&\beta_x y\underline{s}(z)\underline{s}'(z)\lim_{n\to\infty}\dfrac{1}{n}\sum_{k=1}^{p}\mathbf{e}_{k}^\top\mathbf{G}^\top\mathbf{M}\left(\mathbf{I}+\underline{s}(z)\mathbf{R}\mathbf{M}\right)^{-1}\mathbf{G}\mathbf{e}_{k}\cdot\mathbf{e}_{k}^\top\mathbf{G}^\top\mathbf{M}\left(\mathbf{I}+\underline{s}(z)\mathbf{R}\mathbf{M}\right)^{-2}\mathbf{G}\mathbf{e}_{k}\\
+&\displaystyle\lim_{n\to \infty } \dfrac{1}{n} \sum_{k=1}^{p}\dfrac{\partial}{\partial z}\left[\underline{s}(z)\mathbf{e}_{k}^{\top}\mathbf{RM}\left(\mathbf{I}+\underline{s}(z)\mathbf{R}\mathbf{M}\right)^{-1}\mathbf{e}_{k}\cdot \mathbf{e}_{k}^{\top}\left(\mathbf{I}+\underline{s}(z)\mathbf{R}\mathbf{M}\right)^{-1}\mathbf{R}\mathbf{e}_{k}\right]\\
+&\displaystyle\lim_{n\to \infty } \dfrac{1}{n} \sum_{k=1}^{p}\dfrac{\partial}{\partial z}\left[\underline{s}(z)\mathbf{e}_{k}^{\top}\mathbf{RM}\left(\mathbf{I}+\underline{s}(z)\mathbf{R}\mathbf{M}\right)^{-1}\mathbf{M}\mathbf{e}_{k}\cdot \mathbf{e}_{k}^{\top}\mathbf{M}\left(\mathbf{I}+\underline{s}(z)\mathbf{R}\mathbf{M}\right)^{-1}\mathbf{R}\mathbf{e}_{k}\right]\\
+&\displaystyle\lim_{n\to \infty } \dfrac{\beta_{x}}{2n} \sum_{k,i=1}^{p}\dfrac{\partial}{\partial z}\left[\underline{s}(z)g_{ki}^{2}\mathbf{e}_{i}^{\top}\mathbf{G}^{\top}\mathbf{M}\left(\mathbf{I}+\underline{s}(z)\mathbf{R}\mathbf{M}\right)^{-1}\mathbf{e}_{k}\cdot\mathbf{e}_{k}^{\top}\left(\mathbf{I}+\underline{s}(z)\mathbf{R}\mathbf{M}\right)^{-1}\mathbf{Ge}_{i}\right]\\
+&\displaystyle\lim_{n\to \infty } \dfrac{\beta_{x}}{2n} \sum_{k,i=1}^{p}\dfrac{\partial}{\partial z}\left[\underline{s}(z)g_{ki}^{2}\mathbf{e}_{i}^{\top}\mathbf{G}^{\top}\mathbf{M}\left(\mathbf{I}+\underline{s}(z)\mathbf{R}\mathbf{M}\right)^{-1}\mathbf{M}\mathbf{e}_{k}\cdot\mathbf{e}_{k}^{\top}\mathbf{M}\left(\mathbf{I}+\underline{s}(z)\mathbf{R}\mathbf{M}\right)^{-1}\mathbf{Ge}_{i}\right].
\end{align*}

(4): The term V converges weakly to a zero mean Gaussian process in both cases. The process is tight in both cases. The variance-covariance function is
\begin{align*}
&v(z_{1},z_{2}) \\
=&2\left\{\dfrac{\underline{s}'(z_{1})\underline{s}'(z_{2})}{\left[\underline{s}(z_{2})-\underline{s}(z_{1})\right]^{2}}-\dfrac{1}{(z_{1}-z_{2})^{2}}\right\}\\
+&\dfrac{1}{2}\displaystyle\lim_{n \to \infty} \dfrac{1}{n} \displaystyle\sum_{k,\ell=1}^{p}r_{k\ell}^{2}\dfrac{\partial}{\partial z_{1}}\left[\mathbf{e}_{k}^{\top} \mathbf{M}(\mathbf{I}_{p}+\underline{s}(z_{1})\mathbf{RM})^{-1} \mathbf{M}^{-1}\mathbf{e}_{k}\right]\cdot \dfrac{\partial}{\partial z_{2}}\left[\mathbf{e}_{\ell}^{\top}\mathbf{M}(\mathbf{I}_{p}+\underline{s}(z_{2})\mathbf{RM})^{-1}\mathbf{M}^{-1}\mathbf{e}_{\ell}\right]\\
+&\dfrac{1}{2}\displaystyle\lim_{n \to \infty} \dfrac{1}{n} \displaystyle\sum_{k,\ell=1}^{p}r_{k\ell}^{2}\dfrac{\partial}{\partial z_{1}}\left[\mathbf{e}_{k}^{\top}(\mathbf{I}_{p}+\underline{s}(z_{1})\mathbf{RM})^{-1}\mathbf{e}_{k}\right]\cdot \dfrac{\partial}{\partial z_{2}}\left[\mathbf{e}_{\ell}^{\top}(\mathbf{I}_{p}+\underline{s}(z_{2})\mathbf{RM})^{-1}\mathbf{e}_{\ell}\right]\\
+&\dfrac{1}{2}\displaystyle\lim_{n \to \infty} \dfrac{1}{n} \displaystyle\sum_{k,\ell=1}^{p}r_{k\ell}^{2}\dfrac{\partial}{\partial z_{1}}\left[\mathbf{e}_{k}^{\top}\mathbf{M}(\mathbf{I}_{p}+\underline{s}(z_{1})\mathbf{RM})^{-1}\mathbf{M}^{-1}\mathbf{e}_{k}\right]\cdot \dfrac{\partial}{\partial z_{2}}\left[\mathbf{e}_{\ell}^{\top}(\mathbf{I}_{p}+\underline{s}(z_{2})\mathbf{RM})^{-1}\mathbf{e}_{\ell}\right]\\
+&\dfrac{1}{2}\displaystyle\lim_{n \to \infty} \dfrac{1}{n} \displaystyle\sum_{k,\ell=1}^{p}r_{k\ell}^{2}\dfrac{\partial}{\partial z_{2}}\left[\mathbf{e}_{k}^{\top}\mathbf{M}(\mathbf{I}_{p}+\underline{s}(z_{2})\mathbf{RM})^{-1}\mathbf{M}^{-1}\mathbf{e}_{k}\right]\cdot \dfrac{\partial}{\partial z_{1}}\left[\mathbf{e}_{\ell}^{\top}(\mathbf{I}_{p}+\underline{s}(z_{1})\mathbf{RM})^{-1}\mathbf{e}_{\ell}\right]\\
-&\underline{s}'(z_{1})\underline{s}'(z_{2})\displaystyle\lim_{n \to \infty}\dfrac{1}{n}\sum_{k=1}^{p}\left[\mathbf{e}_{k}^{\top}(\mathbf{I}_{p}+\underline{s}(z_{2})\mathbf{RM})^{-2}\mathbf{RMR}\mathbf{e}_{k}\right]\cdot \left[\mathbf{e}_{k}^{\top}(\mathbf{I}_{p}+\underline{s}(z_{1})\mathbf{RM})^{-2}\mathbf{RM}\mathbf{e}_{k}\right]\\
-&\underline{s}'(z_{1})\underline{s}'(z_{2})\displaystyle\lim_{n \to \infty}\dfrac{1}{n}\sum_{k=1}^{p}\left[\mathbf{e}_{k}^{\top}(\mathbf{I}_{p}+\underline{s}(z_{1})\mathbf{RM})^{-2}\mathbf{RMR}\mathbf{e}_{k} \right]\cdot\left[\mathbf{e}_{k}^{\top}(\mathbf{I}_{p}+\underline{s}(z_{2})\mathbf{RM})^{-2}\mathbf{RM}\mathbf{e}_{k}\right]\\
-&\underline{s}'(z_{1})\underline{s}'(z_{2})\displaystyle\lim_{n \to \infty}\dfrac{1}{n}\sum_{k=1}^{p}\left[ \mathbf{e}_{k}^{\top}(\mathbf{I}_{p}+\underline{s}(z_{2})\mathbf{RM})^{-2}\mathbf{RMR}\mathbf{e}_{k}\right]\cdot\left[\mathbf{e}_{k}^{\top}\mathbf{M}(\mathbf{I}_{p}+\underline{s}(z_{1})\mathbf{RM})^{-2}\mathbf{R}\mathbf{e}_{k}\right]\\
-&\underline{s}'(z_{1})\underline{s}'(z_{2})\displaystyle\lim_{n \to \infty}\dfrac{1}{n}\sum_{k=1}^{p}\left[ \mathbf{e}_{k}^{\top}(\mathbf{I}_{p}+\underline{s}(z_{1})\mathbf{RM})^{-2}\mathbf{RMR}\mathbf{e}_{k}\right]\cdot\left[\mathbf{e}_{k}^{\top}\mathbf{M}(\mathbf{I}_{p}+\underline{s}(z_{2})\mathbf{RM})^{-2}\mathbf{R}\mathbf{e}_{k}\right]
\end{align*}
in the elliptical case. While the variance-covariance function in the linear case becomes
\begin{align*}
&v(z_{1},z_{2})\\
=&2\left\{\dfrac{\underline{s}'(z_{1})\underline{s}'(z_{2})}{\left[\underline{s}(z_{2})-\underline{s}(z_{1})\right]^{2}}-\dfrac{1}{(z_{1}-z_{2})^{2}}\right\}\\
+&\lim_{n\to \infty}\dfrac{\beta_x y}{n}\sum_{k=1}^p\dfrac{\partial}{\partial z_{1}}\underline{s}(z_{1})\mathbf{e}_{k}^{\top}\mathbf{G}^\top\mathbf{M}\left(\mathbf{I}_{p}+\underline{s}(z_1)\mathbf{R}\mathbf{M}\right)^{-1}\mathbf{G}\mathbf{e}_k\cdot\dfrac{\partial}{\partial z_2}\underline{s}(z_2)\mathbf{e}_{k}^{\top}\mathbf{G}^\top\mathbf{M}\left(\mathbf{I}_{p}+\underline{s}(z_2)\mathbf{R}\mathbf{M}\right)^{-1}\mathbf{G}\mathbf{e}_k\\
+&\dfrac{1}{4}\displaystyle\lim_{n\to \infty}\dfrac{1}{n}\sum_{k,\ell=1}^{p}\left(\beta_{x}\sum_{i=1}^{p}g_{ki}^{2}g_{\ell i}^{2}+2r_{k\ell}^{2}\right)\cdot\dfrac{\partial}{\partial z_{1}}\mathbf{e}_{k}^{\top}(\mathbf{I}_{p}+\underline{s}(z_{1})\mathbf{RM})^{-1}\mathbf{e}_{k}\cdot\dfrac{\partial}{\partial z_{2}}\mathbf{e}_{\ell}^{\top}(\mathbf{I}_{p}+\underline{s}(z_{2})\mathbf{RM})^{-1}\mathbf{e}_{\ell}\\
+&\dfrac{1}{4}\displaystyle\lim_{n\to \infty}\dfrac{1}{n}\sum_{k,\ell=1}^{p}\left(\beta_{x}\sum_{i=1}^{p}g_{ki}^{2}g_{\ell i}^{2}+2r_{k\ell}^{2}\right)\cdot\dfrac{\partial}{\partial z_{1}}\mathbf{e}_{k}^{\top}\mathbf{M}(\mathbf{I}_{p}+\underline{s}(z_{1})\mathbf{RM})^{-1}\mathbf{M}^{-1}\mathbf{e}_{k}\\
&\quad\quad\quad\quad\quad\quad\quad\cdot\dfrac{\partial}{\partial z_{2}}\mathbf{e}_{\ell}^{\top}\mathbf{M}(\mathbf{I}_{p}+\underline{s}(z_{2})\mathbf{RM})^{-1}\mathbf{M}^{-1}\mathbf{e}_{\ell}\\
+&\dfrac{1}{4}\displaystyle\lim_{n\to \infty}\dfrac{1}{n}\sum_{k,\ell=1}^{p}\left(\beta_{x}\sum_{i=1}^{p}g_{ki}^{2}g_{\ell i}^{2}+2r_{k\ell}^{2}\right)\dfrac{\partial}{\partial z_{1}}\mathbf{e}_{k}^{\top}(\mathbf{I}_{p}+\underline{s}(z_{1})\mathbf{RM})^{-1}\mathbf{e}_{k}\\
&\quad\quad\quad\quad\quad\quad\quad\cdot\dfrac{\partial}{\partial z_{2}}\mathbf{e}_{\ell}^{\top}\mathbf{M}(\mathbf{I}_{p}+\underline{s}(z_{2})\mathbf{RM})^{-1}\mathbf{M}^{-1}\mathbf{e}_{\ell}\\
+&\dfrac{1}{4}\displaystyle\lim_{n\to \infty}\dfrac{1}{n}\sum_{k,\ell=1}^{p}\left(\beta_{x}\sum_{i=1}^{p}g_{ki}^{2}g_{\ell i}^{2}+2r_{k\ell}^{2}\right)\dfrac{\partial}{\partial z_{2}}\mathbf{e}_{k}^{\top}(\mathbf{I}_{p}+\underline{s}(z_{2})\mathbf{RM})^{-1}\mathbf{e}_{k}\\
&\quad\quad\quad\quad\quad\quad\quad\cdot\dfrac{\partial}{\partial z_{1}}\mathbf{e}_{\ell}^{\top}\mathbf{M}(\mathbf{I}_{p}+\underline{s}(z_{1})\mathbf{RM})^{-1}\mathbf{M}^{-1}\mathbf{e}_{\ell}\\
-&\displaystyle\lim_{n\to \infty}\dfrac{1}{n}\sum_{k=1}^{p}\dfrac{\partial}{\partial z_{1}}\mathbf{e}_{k}^{\top}(\mathbf{I}_{p}+\underline{s}(z_{1})\mathbf{RM})^{-1}\mathbf{e}_{k}\cdot \dfrac{\partial}{\partial z_{2}}\mathbf{e}_{k}^{\top}(\mathbf{I}_{p}+\underline{s}(z_{2})\mathbf{RM})^{-1}\mathbf{R}\mathbf{e}_{k}\\
-&\displaystyle\lim_{n\to \infty}\dfrac{1}{n}\sum_{k=1}^{p}\dfrac{\partial}{\partial z_{2}}\mathbf{e}_{k}^{\top}(\mathbf{I}_{p}+\underline{s}(z_{2})\mathbf{RM})^{-1}\mathbf{e}_{k}\cdot \dfrac{\partial}{\partial z_{1}}\mathbf{e}_{k}^{\top}(\mathbf{I}_{p}+\underline{s}(z_{1})\mathbf{RM})^{-1}\mathbf{R}\mathbf{e}_{k}\\
+&\displaystyle\lim_{n\to \infty}\dfrac{\beta_{x}}{2n}\sum_{k,\ell=1}^{p}g_{k\ell}^{2}\dfrac{\partial}{\partial z_{1}}\mathbf{e}_{k}^{\top}(\mathbf{I}_{p}+\underline{s}(z_{1})\mathbf{RM})^{-1}\mathbf{e}_{k}\cdot \dfrac{\partial}{\partial z_{2}}\underline{s}(z_{2})\mathbf{e}_{\ell}^{\top}\mathbf{G}^{\top}(\mathbf{I}_{p}+\underline{s}(z_{2})\mathbf{RM})^{-1}\mathbf{G}\mathbf{e}_{\ell}\\
+&\displaystyle\lim_{n\to \infty}\dfrac{\beta_{x}}{2n}\sum_{k,\ell=1}^{p}g_{k\ell}^{2}\dfrac{\partial}{\partial z_{2}}\mathbf{e}_{k}^{\top}(\mathbf{I}_{p}+\underline{s}(z_{2})\mathbf{RM})^{-1}\mathbf{e}_{k}\cdot \dfrac{\partial}{\partial z_{1}}\underline{s}(z_{1})\mathbf{e}_{\ell}^{\top}\mathbf{G}^{\top}(\mathbf{I}_{p}+\underline{s}(z_{1})\mathbf{RM})^{-1}\mathbf{G}\mathbf{e}_{\ell}\\
-&\displaystyle\lim_{n\to \infty}\dfrac{1}{n}\sum_{k=1}^{p}\dfrac{\partial}{\partial z_{1}}\mathbf{e}_{k}^{\top}\mathbf{M}(\mathbf{I}_{p}+\underline{s}(z_{1})\mathbf{RM})^{-1}\mathbf{M}^{-1}\mathbf{e}_{k}\cdot \dfrac{\partial}{\partial z_{2}}\mathbf{e}_{k}^{\top}(\mathbf{I}_{p}+\underline{s}(z_{2})\mathbf{RM})^{-1}\mathbf{R}\mathbf{e}_{k}\\
-&\displaystyle\lim_{n\to \infty}\dfrac{1}{n}\sum_{k=1}^{p}\dfrac{\partial}{\partial z_{2}}\mathbf{e}_{k}^{\top}\mathbf{M}(\mathbf{I}_{p}+\underline{s}(z_{2})\mathbf{RM})^{-1}\mathbf{M}^{-1}\mathbf{e}_{k}\cdot \dfrac{\partial}{\partial z_{1}}\mathbf{e}_{k}^{\top}(\mathbf{I}_{p}+\underline{s}(z_{1})\mathbf{RM})^{-1}\mathbf{R}\mathbf{e}_{k}\\
+&\displaystyle\lim_{n\to \infty}\dfrac{\beta_{x}}{2n}\sum_{k,\ell=1}^{p}g_{k\ell}^{2}\dfrac{\partial}{\partial z_{1}}\mathbf{e}_{k}^{\top}\mathbf{M}(\mathbf{I}_{p}+\underline{s}(z_{1})\mathbf{RM})^{-1}\mathbf{M}^{-1}\mathbf{e}_{k}\cdot \dfrac{\partial}{\partial z_{2}}\underline{s}(z_{2})\mathbf{e}_{\ell}^{\top}\mathbf{G}^{\top}(\mathbf{I}_{p}+\underline{s}(z_{2})\mathbf{RM})^{-1}\mathbf{G}\mathbf{e}_{\ell}\\
+&\displaystyle\lim_{n\to \infty}\dfrac{\beta_{x}}{2n}\sum_{k,\ell=1}^{p}g_{k\ell}^{2}\dfrac{\partial}{\partial z_{2}}\mathbf{e}_{k}^{\top}\mathbf{M}(\mathbf{I}_{p}+\underline{s}(z_{2})\mathbf{RM})^{-1}\mathbf{M}^{-1}\mathbf{e}_{k}\cdot \dfrac{\partial}{\partial z_{1}}\underline{s}(z_{1})\mathbf{e}_{\ell}^{\top}\mathbf{G}^{\top}(\mathbf{I}_{p}+\underline{s}(z_{1})\mathbf{RM})^{-1}\mathbf{G}\mathbf{e}_{\ell}.
\end{align*}

Combining (1)-(4) and (\ref{eqA}), we conclude that the random vector
\begin{equation}
\left(\displaystyle\sum_{i=1}^{p}g_{1}(\hat{\lambda}_{i})-p\displaystyle\int g_{1}(x)\mathrm{d}F^{y_{n},H_{n}}(x),\ldots,\displaystyle\sum_{i=1}^{p}g_{K}(\hat{\lambda}_{i})-p\displaystyle\int g_{K}(x)\mathrm{d}F^{y_{n},H_{n}}(x)\right)
\end{equation}
converges to a K-dimensional normal random vector $(X_{g_{1}},\ldots, X_{g_{K}})$ in both cases. The mean function is
$$
\textup{E}X_{g_{\ell}}=-\dfrac{1}{2\pi i}\oint g_{\ell}(z)(\textup{E}_{0}(z)+\textup{E}_{1}(z))\mathrm{d}z,\quad 1\leq \ell\leq K.
$$
The variance-covariance function is
$$
\textup{Cov}(X_{g_{\ell_{1}}},X_{g_{\ell_{2}}})=-\dfrac{1}{4\pi ^{2}}\oint \hspace{-0.2cm}\oint g_{\ell_{1}}(z_{1})g_{\ell_{2}}(z_{2})v(z_{1},z_{2})\mathrm{d}z_{1}\mathrm{d}z_{2},\quad 1\leq \ell_{1} \leq \ell_{2} \leq K.
$$
This completes the proof of Theorems \ref{th3.1} and \ref{th4.1}.

\subsection{Proof of Theorem \ref{th3.2} and Example \ref{ex3.1}}
\subsubsection{Proof of Theorem \ref{th3.2}}
We take the same approach as in \cite{yin2023central}. Let $\underline{s}(z)=\frac{1}{1+\sqrt{y}\xi},$ that means we also regard $z$ as a function of $\xi$. Then Theorem \ref{th3.2}
 follows.

\subsubsection{Proof of Example \ref{ex3.1}}
First, we know that the moments of the standard M-P distribution with index y take the values
$$
m_k(y)=\sum_{r=0}^{k-1} \dfrac{1}{r+1}\binom{k}{r}\binom{k-1}{r} y^k
$$
see Lemma 3.1 in \cite{bai2010spectral}. From this, we can calculate the centering terms
$$
\int g_1(x) f^{y_{n-1}}(x) \mathrm{d} x=m_1\left(y_{n-1}\right)=1, \quad \int g_2(x) f^{y_{n-1}}(x) \mathrm{d} x=m_2\left(y_{n-1}\right)=1+y_{n-1}.
$$
Then, we note that the following equation holds
$$\dfrac{\xi^2-1}{\xi^3(\xi+\sqrt{y})}=\dfrac{y-1}{y\sqrt{y}\xi}+\dfrac{1}{y\xi^2}-\dfrac{1}{\sqrt{y}\xi^3}+\dfrac{1-y}{y\sqrt{y}(\xi+\sqrt{y})}.$$
Denote
\begin{align*}
    \mu_1(g)&=\displaystyle\lim_{r\to 1^{+}}\dfrac{1}{2\pi i}\oint_{|\xi|=1}g(|1+\sqrt{y}\xi|^2)\left(\dfrac{\xi}{\xi^2-r^{-2}}-\dfrac{1}{\xi}\right) \mathrm{d}\xi,\\
    \mu_2(g)&=\dfrac{1}{2\pi i}\oint_{|\xi|=1}g(|1+\sqrt{y}\xi|^2)\dfrac{1}{\xi^3} \mathrm{d}\xi,\\
     \mu_3(g)&=\dfrac{1}{2\pi i}\oint_{|\xi|=1}g(|1+\sqrt{y}\xi|^2)\dfrac{1}{\xi^2} \mathrm{d}\xi,\\
      \mu_4(g)&=\dfrac{1}{2\pi i}\oint_{|\xi|=1}g(|1+\sqrt{y}\xi|^2)\dfrac{1}{\xi} \mathrm{d}\xi,\\
       \mu_5(g)&=\dfrac{1}{2\pi i}\oint_{|\xi|=1}g(|1+\sqrt{y}\xi|^2)\dfrac{1}{\xi+\sqrt{y}} \mathrm{d}\xi.\\
\end{align*}
By the Residue theorem and \cite{yin2022spectral}, we can calculate that
\begin{align*}
    &\mu_1(g_1)=0,~\mu_2(g_1)=0,~\mu_3(g_1)=\sqrt{y}, ~\mu_4(g_1)=1+y,~\mu_5(g_1)=1,\\
    &\mu_1(g_2)=y,~\mu_2(g_2)=y,~\mu_3(g_2)=2\sqrt{y}(1+y), ~\mu_4(g_2)=y^2+4y+1,\\
    &\mu_5(g_2)=2y+1.
\end{align*}
The results of this example follows.

\newpage

\bibliographystyle{siam}
\bibliography{library}

\end{document}